\documentclass[10pt,reqno]{amsart}
  \usepackage{geometry}
\usepackage[T1, T2A]{fontenc}
\usepackage[utf8]{inputenc} 
\usepackage{graphicx}
\usepackage{amssymb,latexsym}
\usepackage{subfigure}
\usepackage{tikz}

\numberwithin{equation}{section}

\def\eps{{\varepsilon}}

\def\C{{\mathbb C}}
\def\D{{\mathbb D}}

\def\N{{\mathbb N}}

\def\R{{\mathbb R}}

\def\Z{{\mathbb Z}}
\def\CP{{\mathbb C\mathbb P}}

\def\ov{\overline}

\usepackage{cleveref}
\theoremstyle{plain}
\newtheorem{lemma}{Lemma}[section]

\theoremstyle{definition}
\newtheorem{example}[lemma]{Example}
\newtheorem{remark}[lemma]{Remark}
\theoremstyle{plain}

\newtheorem{theorem}[lemma]{Theorem}

\theoremstyle{definition}

\theoremstyle{remark}

\title[The equivalence problem for $1$-resonance]{The equivalence problem in analytic dynamics for $1$-resonance\footnote{The author is supported by NSERC in Canada.}}

\author[C. Rousseau]{Christiane Rousseau}
\address{Christiane Rousseau, D\'epartement de
math\'ematiques et de statistique, Universit\'e de Montr\'eal, C.P. 6128,
Succursale Centre-ville, Montr\'eal (Qc), H3C 3J7, Canada.}
\email{christiane.rousseau@umontreal.ca}

\usepackage{xcolor}
\usepackage{xspace}


\begin{document}

\date{\today}

\begin{abstract} 
When are two germs of analytic systems conjugate or orbitally equivalent under an analytic change of coordinates in the neighborhood of a singular point?
A way to answer is to use normal forms. But there are large classes of dynamical systems for which the change of coordinates to a normal form diverges. In this paper we discuss the case of singularities for which the normalizing transformation is $k$-summable, thus allowing to provide moduli spaces. We explain the common geometric features of these singularities, and show that the study of their unfoldings allows understanding both the singularities themselves, and the geometric obstructions to convergence of the normalizing transformations. We also present some moduli spaces for generic $k$-parameter families unfolding such singularities. 
\end{abstract}

\maketitle

\section{Introduction}

Singularities of dynamical systems organize the dynamics, thus explaining why their study is so important. In generic situations, usually only simple singularities occur. But dynamical systems often depend on parameters and, the more parameters in the system, the more complex the singularities. 
In \cite{A}, Vladimir Arnold explains that the singularities of codimension $\leq k$ are unavoidable in $k$-parameter families of dynamical systems, while we can get rid of singularities of codimension greater than $k$ by slightly perturbing the family.

In this paper, we discuss singularities of  analytic dynamical systems. A fundamental problem is the equivalence problem: \emph{When are two germs of analytic dynamical systems 
equivalent in the neighborhood of a singularity under an analytic
change of coordinates?}
One way of solving the equivalence problem is to use normal forms, the simplest case being when there exists a linearizing change of coordinates. 
It is natural to look for a normalizing change of coordinate as a power series and then to study its convergence. But the normalizing changes of coordinates only converge for the simplest singularities, for instance for a fixed point of a $1$-dimensional diffeomorphism with a multiplier of norm different from $1$, or for a node of a planar vector field. Even in the case of a saddle of a planar vector field (which is a hyperbolic singularity) the change to normal form generically diverges as soon as the ratio of eigenvalues is not a diophantian irrational number.  

The question we are interested in this paper is: 
\smallskip

\centerline{\bf \Large Why?} 
\smallskip

\noindent Why is it so often the case that the change of coordinates to normal form diverges? 

Let us start by discussing by discussing the two examples mentioned above where we have convergence to the normal form. 
The first one is that of a fixed point of a $1$-dimensional diffeomorphism defined on a neighborhood $V$ of the origin with a multiplier $\lambda$ such that  $|\lambda|\neq1$. The orbit space is the quotient $V/f \cup\{0\}$, which is the union of the origin with a torus of modulus $\frac{\log \lambda}{2\pi i}$. The torus has a unique complex structure, which is independent of the size of $U$ and of the special form of the diffeomorphism. Hence, all such $1$-dimensional  germs of diffeomorphims are conjugate. 

The case of the orbital normal form of a planar vector field $X$ in the neighborhood of a node with eigenvalues $1$ and $\lambda\in [1, +\infty)$ is similar. In that case, we have convergence to the normal form $X'$, which is linear as soon as $\lambda\notin\N$: $X'= x \frac{\partial}{\partial x} + \lambda y \frac{\partial}{\partial y}$, and the solutions are given by $y = Cx^\lambda$, $C\in \C$, together with $x=0$, which corresponds to the limit case $C=\infty$. The complex space of leaves of the underlying foliation of the linear node is a quotient $\CP^1/L_B$, where $L_B$ is the linear map $C\mapsto BC$ for $B=\exp(2\pi i \lambda)$. This very rigid object with a unique complex structure is also the space of leaves of any node with eigenvalues $1$ and $\lambda$ when $\lambda\notin \N$. Hence, this forces the convergence of the orbital linearization. When $\lambda=n\in\N$, the node is resonant and the orbital normal form is given by $X'= x \frac{\partial}{\partial x} + (\lambda y +Ay^n)\frac{\partial}{\partial y}$, for $A\in\{0,1\}$. When $A=0$, the space of leaves is the whole of $\CP^1$, while for $A=1$ the space of leaves is $\CP^1/T_D$ for the translation $T_D(C) = C+D$ with $D=2\pi i$; indeed, the leaves have the form $y= x^n\log x + Cx^n$, $C\in \C$, together with $x=0$. Again, the space of leaves is a very rigid object. Even if we had started with a node of a real planar vector field, we see that the geometric explanation of the convergence of the normalizing transformation is seen when extending the system to $(\C^2,0)$. 

In this paper we will discuss several singularities, $1$-resonant singularities, which share the common properties that the infinite number of resonances in the normal form all come from one rational relation between the eigenvalues, and that the formal normal form contains only a finite number of parameters. These singularities have been studied in the literature by \'Ecalle \cite{E}, Voronin\cite{V}, Martinet-Ramis (\cite{MR1}, \cite{MR2}), etc.,  and it has been shown that the change of coordinates to the normal form is generically divergent and $k$-summable. (For the reader not familiar with $k$-summability, a $k$-summable series $\sum_{n\in \N} a_nz^n$ is such that $|a_n| \leq \frac{M (n!)^{1/k}}{r^n}$ for some positive $M,r$,  and we can find unique \lq\lq sums\rq\rq\ asymptotic to the power series on sectors of opening larger than $\pi/k$ covering a neighborhood of the origin. Since we will not work with $k$-summability we do not give a more precise definition.) The divergence of the normalizing transformation comes from the fact that the \lq\lq geometry of the system\rq\rq\   is more complicated than the geometry of the formal normal form: the normal form is too poor to encode all the complex phenomena that can occur. 

A common feature of the $1$-resonant singularities is that they correspond, in the codimension $k$ case, to the coallescence of $k+1$ special \lq\lq objects\rq\rq, which could be fixed points, periodic orbits, special leaves, etc. Because we are studying limit situations it is natural to unfold the systems: when considering a singularity of codimension $k$, its full richness can only be uncovered by studying a generic $k$-parameter unfolding. A large program has been started around 2000 to systematically solve the equivalence problem for generic $k$-parameter unfoldings of  1-resonant singularities of dynamical systems of codimension $k$, to which the author has contributed with her collaborators and students:  Arriagada-Silva, Christopher, Hurtubise,  Klime\v{s}, Lambert, Marde\v{s}i\'c, Roussarie and Teyssier (see for instance \cite{MRR}, \cite{Ro1}, \cite{RT1}, \cite{Ro3},  \cite{AR}, \cite{LR},  \cite{CR}, \cite{Ki},\cite{HLR} and \cite{HR}), as well as Rib\'on (\cite{Ri1} and \cite{Ri2}). The idea of unfolding to study the coallescence was not new. In the case of $1$-resonant systems, it has been proposed by several mathematicians including Arnold and Bolibruch, and partially studied by Martinet \cite{Ma}, Ramis \cite{Ra}, Duval \cite{Du} and Glutsyuk (\cite{Gl1} and \cite{Gl2}). The underlying idea of these pioneering works was to unfold in sectors in parameter space where the singularities are hyperbolic. In the neighborhood of the singularities the system can be analytically normalized in an almost rigid way. The comparison of these local normalizations, whenever the domains of normalizations do intersect, is a measure of the obstructions to a global normalization of the system. If these obstructions have a nontrivial limit at the confluence, then we have divergence of the normalizing transformation at the confluence. The limit of this approach in that it does not allow to treat the parameter values for which, either the local normalizations do not exist, or they exist, but their domains do not intersect. For instance, when unfolding a saddle-node, the saddle has no convergent normal form depending on parameters. 
A breakthrough in this program occurred with the vision of Adrien Douady. The thesis of his student, Lavaurs (see \cite{La}), permitted to treat the unfolding of the parabolic point (double fixed point of a 1-dimensional diffeomorphism) in a sector complementary to the one studied by Martinet. The study by Douady, Estrada and Sentenac of the dynamics of complex polynomial vector fields on $\CP^1$ (see \cite{DES}) provided the geometric tools for the equivalence problem in  the general case. 

The paper is organized as follows. In Section~\ref{sec:Examples}, we describe seven examples of $1$-resonant singularities for which the normalizing transformation generically diverges and is $k$-summable. In Section~\ref{sec:sn}, we explain in detail one geometric obstruction to the existence of an analytic normalizing transformation in the case of the saddle-node. In Section~\ref{sec:common}, we present the general features common to all seven examples, and we revisit the examples in more detail. The particular case of the nonresonant irregular singular point of Poincar\'e rank $k$ of a linear differential system is discussed in detail in Section~\ref{sec:Example7}, and we present the kind of results that can be obtained in solving the equivalence problem for generic unfoldings. We end up with a short section of perspectives.

\section{Examples}\label{sec:Examples}
We present here the seven examples that we will follow along the paper. In all these examples the change of coordinate to normal
form is generically $1$-summable.

\subsection{Example 1:}\label{sec:ex1} 
\begin{figure}\begin{center}
\subfigure[$k=1$]{\includegraphics[height=4.5cm]{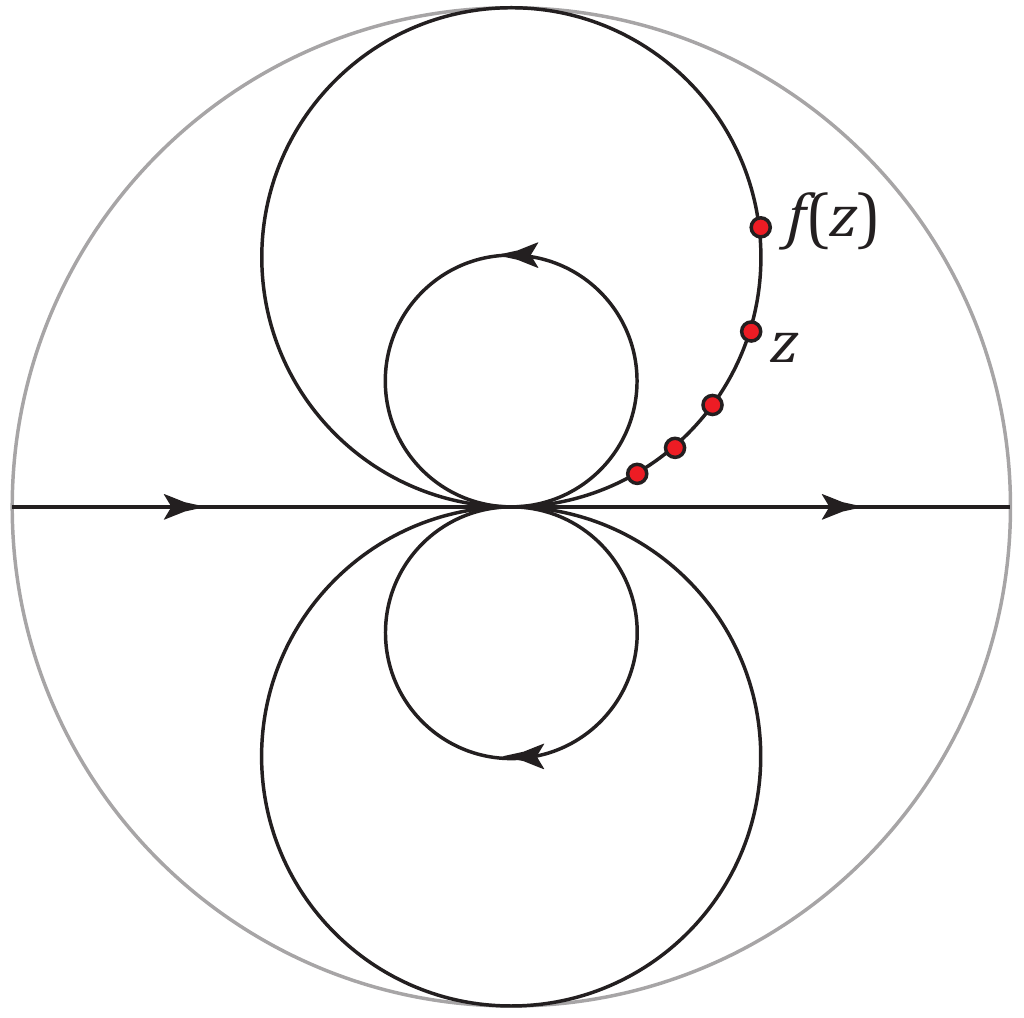}}\qquad\subfigure[$k=4$]{\includegraphics[height=4.5cm]{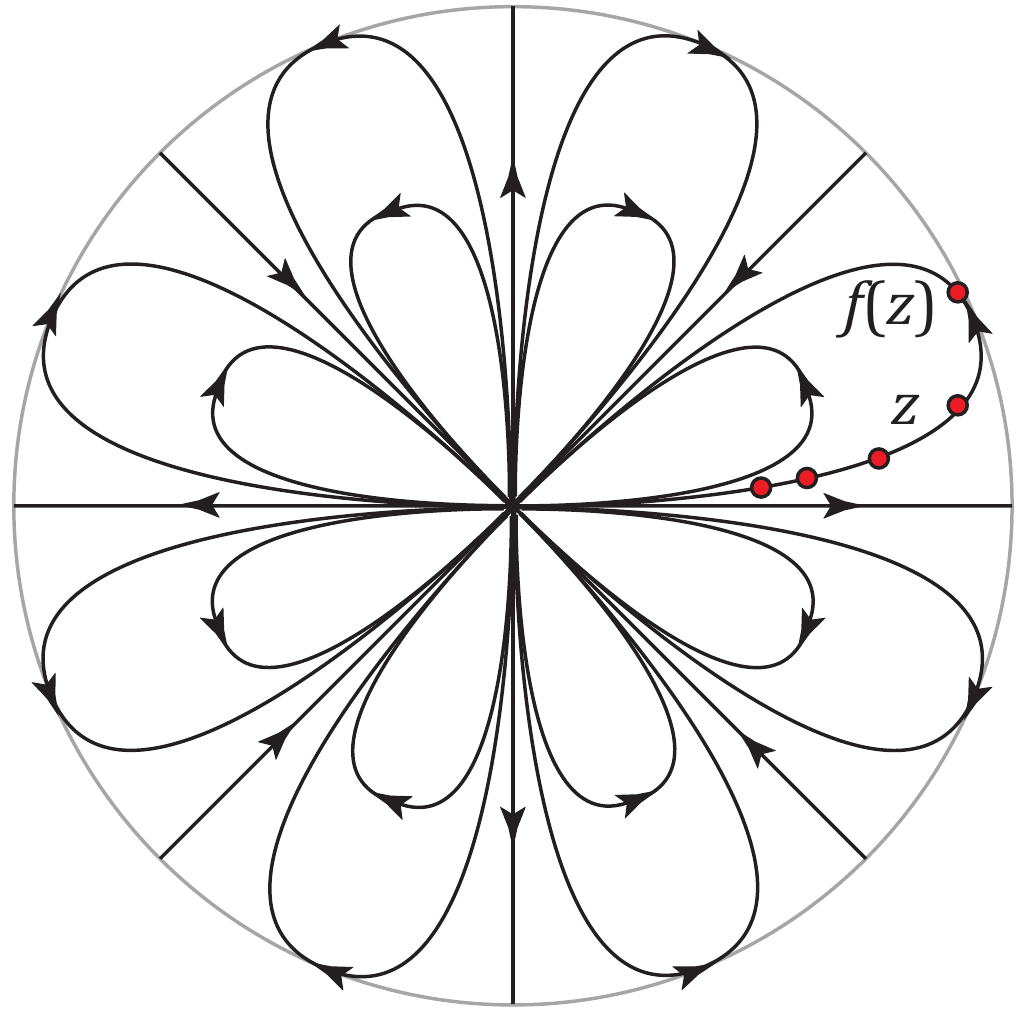}} \end{center}\caption{Parabolic points of codimension $1$ and $4$.}\label{fig:parabolic}\end{figure}

A germ of analytic diffeomorphism $f:
(\C,0)\rightarrow (\C,0)$  with a parabolic point of codimension $k$ (see Figure~\ref{fig:parabolic}) such that  \begin{equation}f(z) =
z+z^{k+1}+\left(\frac{k+1}2-a\right)z^{2k+1}+ o(z^{2k+1}).\label{par_germ}\end{equation} 
Such a point is a multiple fixed point of multiplicity $k+1$, i.e. the coallescence of $k+1$ simple fixed points. 
The formal normal form is the time-one map of the vector field $$\dot z= \frac{z^{k+1}}{1+az^k}.$$

As we will see below, this example is underlying many others. Hence we give more detail and explain how to solve the classification problem for germs of the form \eqref{par_germ}. There are two formal invariants, the codimension, $k$, and the parameter $a$, called in the literature \emph{formal invariant} or \emph{r\'esidu it\'eratif}. The geometric meaning of $a$ will be discussed later. 

\begin{figure}\begin{center}
\subfigure[Left domain]{\includegraphics[width=3.5cm]{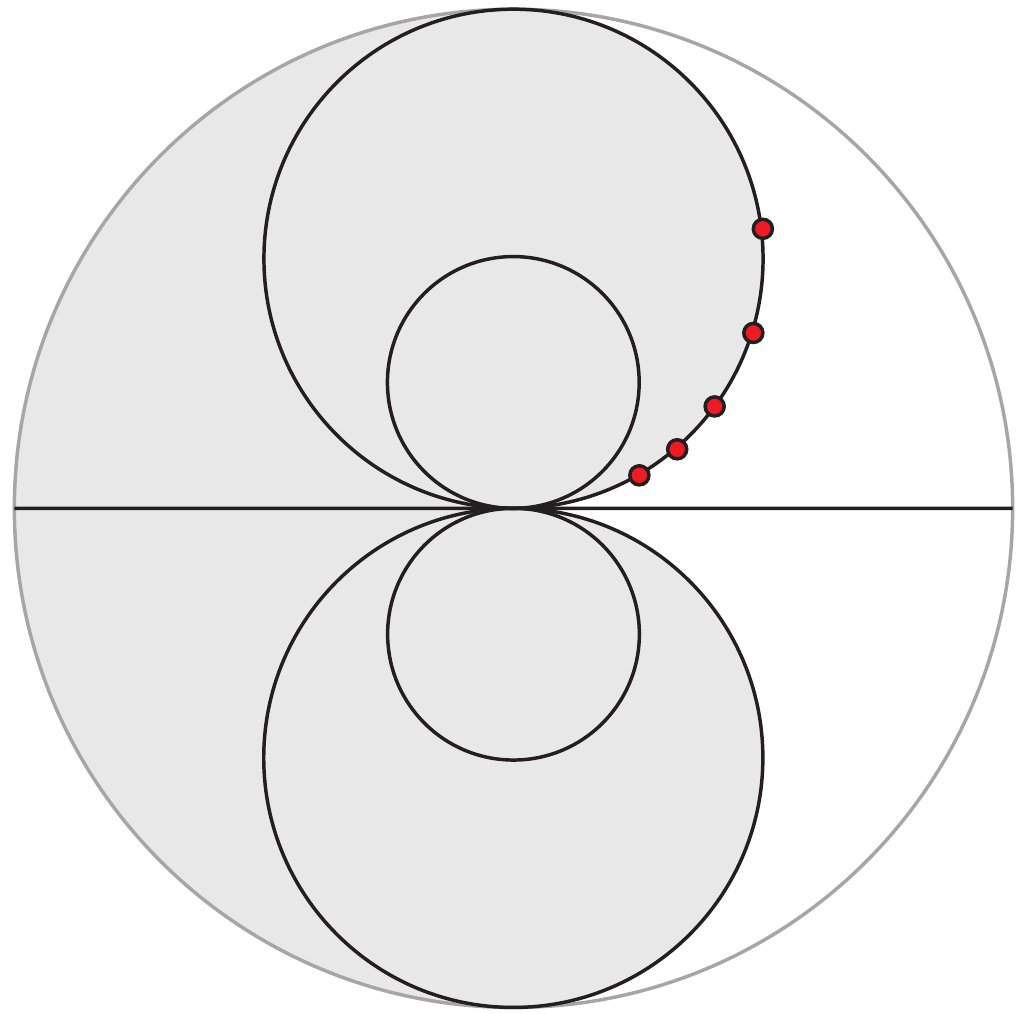}} \quad\subfigure[Right domain]{\includegraphics[width=3.5cm]{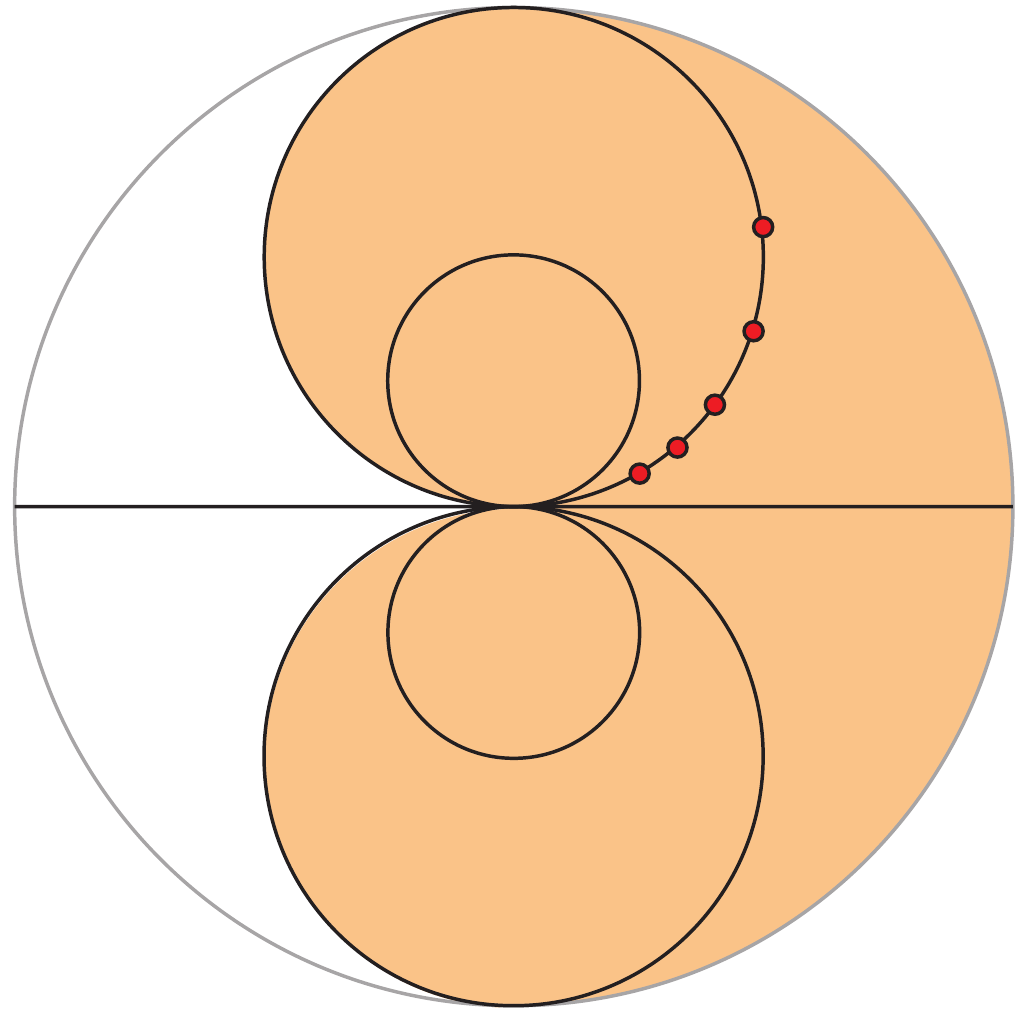}}\quad\subfigure[The intersection]{\includegraphics[width=3.5cm]{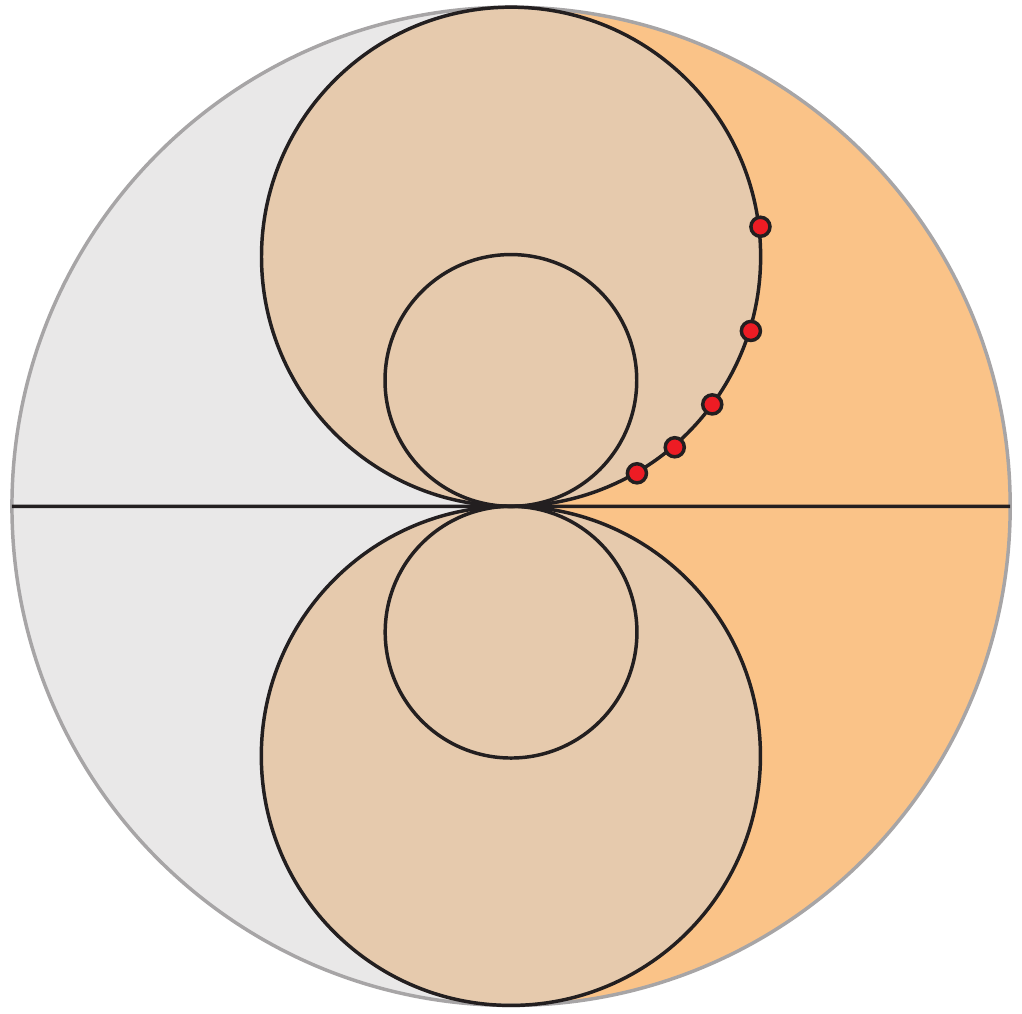}}\caption{The two domains on which almost unique normalizing changes of coordinates exist.}\label{fig:domains_parabolic}\end{center}
\end{figure}

Let us now specialize to the case $k=1$. It is possible to find almost unique normalizing changes of coordinates on two  domains as in Figure~\ref{fig:domains_parabolic}. Because of this rigidity, the mismatch between  the changes of coordinates on the intersection of the domains is the analytic part of the invariant. Note that the intersection has two parts (Figure~\ref{fig:domains_parabolic}(c)). To give this invariant in practice, we will describe the space of orbits. 

\begin{figure}\begin{center}\includegraphics[height=5.5cm]{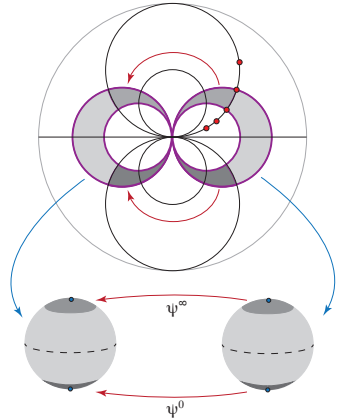} \end{center}\caption{The orbit space of a germ of codimension $1$ parabolic diffeomorphim.}\label{fig:parabolic_orbit}\end{figure}

 We would like to have one representative of each orbit $\{f^{\circ n}(z)\}_{n\in\Z}$, where we limit ourselves to the iterates that exist and lie in a disk $\D_r$. For that purpose we take a curve $\ell$ transversal to the flow lines of $\dot z= \frac{z^2}{1+az}$, and its image $f(\ell)$ (see Figure~\ref{fig:parabolic_orbit}). These two lines form the boundary of a crescent. If we identify the two boundaries, this crescent is conformally equivalent to a sphere with two distinguished points, $0$ and $\infty$, corresponding to the parabolic point. Any point of the sphere represents an orbit, but there are some points of the disk whose orbit does not intersect the crescent. Hence, we repeat the process with a second curve $\ell'$ and its image $f(\ell')$. Again, the crescent bounded by these two curves is conformally equivalent to a sphere with two distinguished points $0$ and $\infty$. But, now there are some orbits which are represented on each of the two spheres. Hence we must identify them in the neighborhoods of $0$ and $\infty$. The identification is done by two germs of analytic diffeomorphisms $\psi^0$ and $\psi^\infty$, often called \emph{horn maps} in the literature. These diffeomorphisms are unique up to the choice of coordinates on the spheres fixing $0$ and $\infty$, i.e. linear changes of coordinates: such changes of coordinates induce an equivalence relation on pairs $(\psi^0,\psi^\infty)$. The equivalence class of a pair of germs $(\psi^0,\psi^\infty)$, noted $[(\psi^0,\psi^\infty)]$ is the analytic part of the modulus.

\begin{theorem}[ \'Ecalle \cite{E} and Voronin \cite{V}] Two parabolic germs of the form \eqref{par_germ} are analytically conjugate, if and only if 
they have the same modulus $\left(k,a,[(\psi^0,\psi^\infty)]\right)$. 
Moreover, any modulus $\left(k,a,[(\psi^0,\psi^\infty)]\right)$ is realizable provided $(\psi^0)'(0)(\psi^\infty)'(\infty)= e^{4\pi^2a}$.\end{theorem}

This modulus space is huge and a parabolic germ is analytically conjugate to its normal form if and only if $\psi^0$ and $\psi^\infty$ are both linear. Hence, convergence to the normal form is an infinite codimension case, which means that there exist a lot of obstructions to convergence! 

\subsection{Example 2:}  
\begin{figure}\begin{center}
\includegraphics[ height=4.5cm]{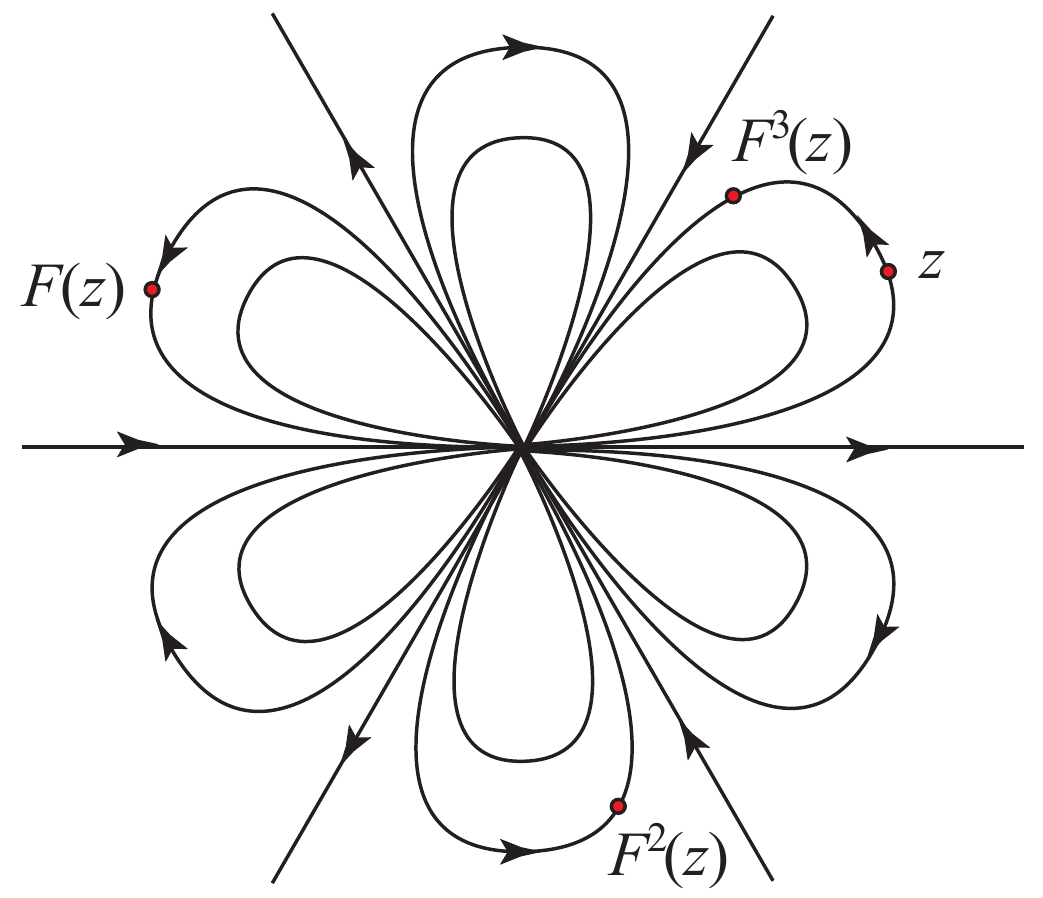}\end{center}\caption{A resonant fixed point of codimension $1$.}\label{fig:periodic}\end{figure}

A germ of analytic
diffeomorphism with a resonant fixed point of codimension $k$: $f: (\C,0)\rightarrow (\C,0)$ such
that (see Figure~\ref{fig:periodic})
$$f(z) = \exp\left(\frac{2\pi i p}{q} \right)\left(z+\frac1{kq}z^{kq+1}+Az^{2kq+1}+o(z^{2kq+1})\right).$$ 
Here the fixed point is of multiplicity $1$, but we see that the origin is a fixed point of multiplicity $kq+1$ of the $q$-th iterate $f^{\circ q}$ of $f$:
$$f^{\circ q}(z) = z+z^{kq+1}+Bz^{2kq+1}+o(z^{2kq+1}),$$ 
where $B=qA+\frac{(kq+1)(q-1)}{2q}$. This is because the origin is the coallescence of one fixed point with $k$ periodic orbits of period $q$. 
The formal normal form for $f^{\circ q}$ is the time-one map of the vector field $\dot z= \frac{z^{kq+1}}{1+az^{kq}}$, where  $A=-\frac{a}{q}+\frac{kq+1}{2q^2}$, and the formal normal form for $f$ is the time-$1/q$ map of the same vector field followed by the rotation of  angle $\frac{2\pi p}{q}$. 

Again, in this example we have the coallescence of $k+1$ objects, which are finite orbits.

\subsection{Example 3:} 
\begin{figure}\begin{center} \includegraphics[width=5cm]{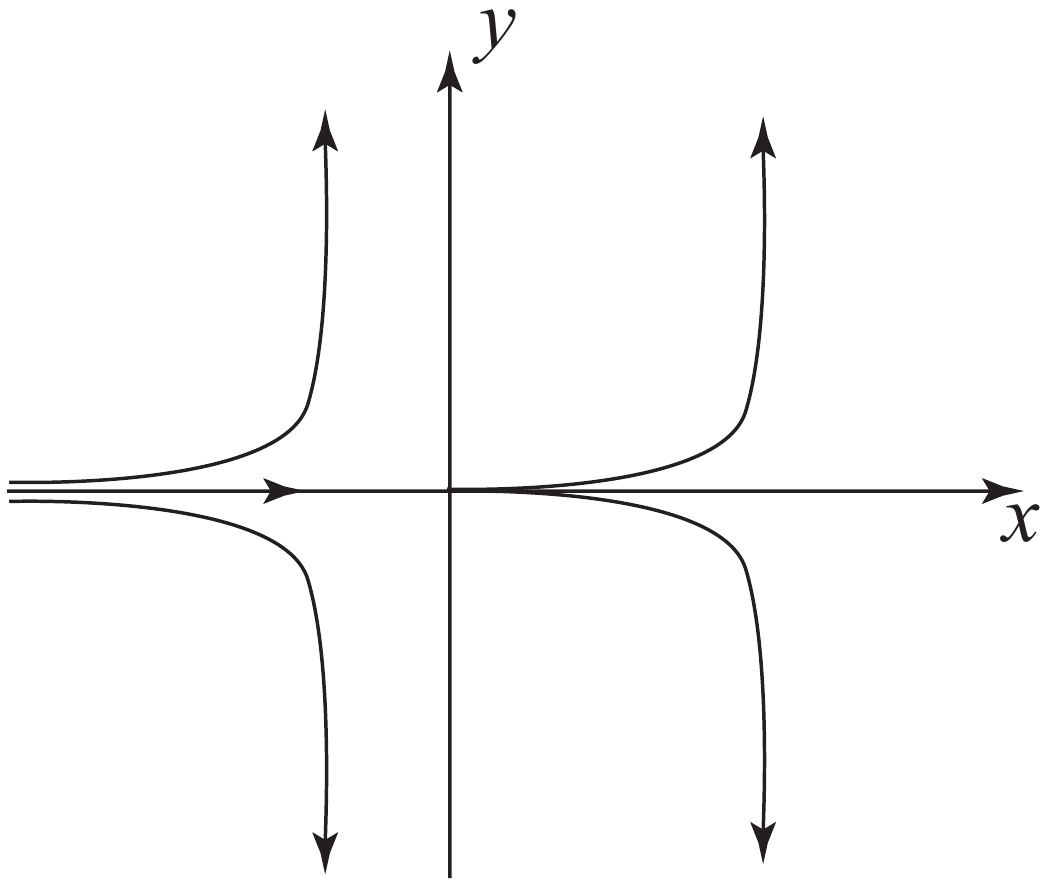}\caption{A planar saddle-node.} \label{fig:saddle-node}\end{center}\end{figure}

\begin{figure}\begin{center} \includegraphics[width=2.5cm]{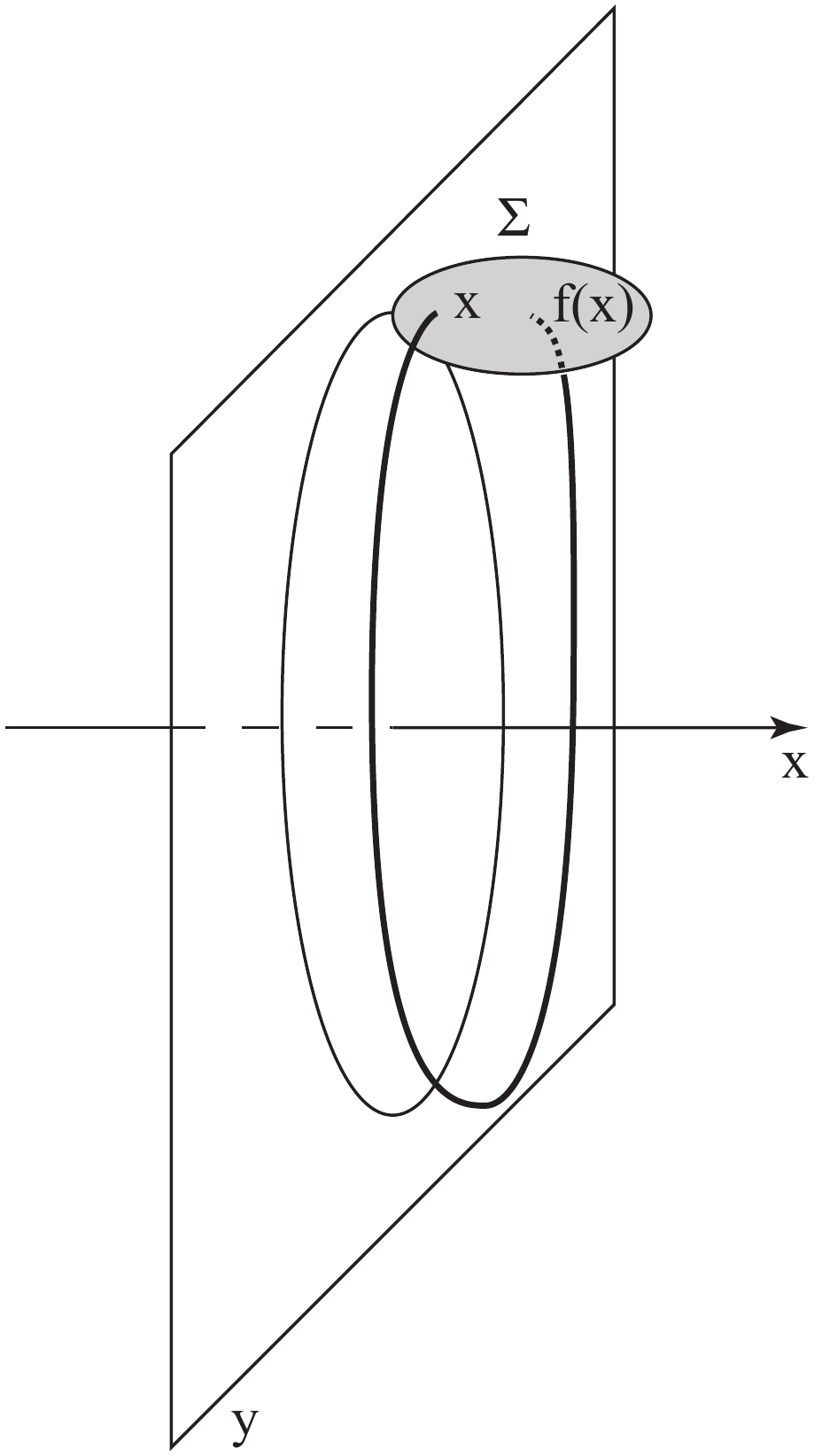}\caption{The holonomy map of a saddle-node}\label{fig:holonomy}\end{center}\end{figure}

A germ of codimension $k$ saddle-node
of a planar vector field (see Figure~\ref{fig:saddle-node})
 with orbital normal form 
$$\begin{cases}
\dot x = x^{k+1},\\
\dot y = y(1+ ax^k).\end{cases}$$
A saddle-node is obviously a singular point of multitplicity $k+1$. There exists a closer relationship with the parabolic point of codimension $k$, which comes from the \emph{holonomy map} of the strong separatrix.
Indeed, a planar vector field defines a local singular foliation in $(\C^2,0)$.  The strong manifold is analytic and can be transformed through a change of coordinate to the coordinate axis $x=0$. The holonomy map is defined on a small neighborhood $\Sigma$ of $x=0$ inside a section $\{y=C\}$ of the strong manifold $x=0$ (see Figure~\ref{fig:holonomy}). 
Consider a closed loop $\gamma=\{(0,C^{ei\theta})\mid \theta\in [0,2\pi]\}$ through $(0,C)$ surrounding the origin in the strong manifold. Let $(x,C)\in \Sigma$ and let $\gamma_x=\{(g(x,\theta),C^{ei\theta})\mid \theta\in [0,2\pi]\}$ be the lifting of $\gamma$ inside the leaf through $(x,C)$. Then the \emph{holonomy of the strong separatrix} is defined as the map 
$$f(x) = g(x,2\pi).$$
This map is analytic and has a parabolic point at $x=0$ of the same codimension $k$ as the saddle-node.
Martinet and Ramis \cite{MR1} proved that two germs of planar vector fields with a saddle-node are orbitally equivalent if and
only if the holonomies of their  strong separatrices are conjugate.

\subsection{Example 4:}
\begin{figure}\begin{center} \includegraphics[width=5cm]{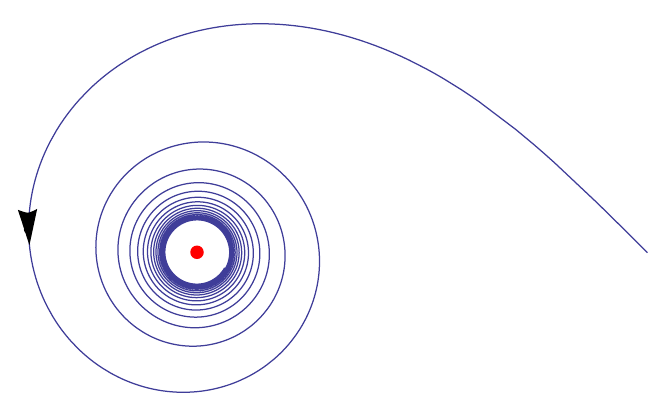}\caption{A weak focus when $k=1$.}\label{fig:weak-focus}\end{center} \end{figure}

A weak focus of order $k$ (see Figure~\ref{fig:weak-focus}), i.e. a singular point of a real analytic planar vector field with pure imaginary eigenvalues, and which is not a center. The orbital normal form is given (up to time reversal) by $$\dot
z = iz -z^{k+1}\overline{z}^k + az^{2k+1}\overline{z}^{2k},$$ with $a\in\R$, which can be rewritten in polar coordinates
$$\begin{cases} \dot r= -r^{2k+1}+ar^{4k+1},\\
\dot \theta = 1.\end{cases}$$
Here again there is a natural $1$-dimensional map associated to the vector field, namely the half Poincar\'e return map defined on a real analytic transversal section passing through the origin. If we parametrize the section by $\zeta$, then this map has the form $$P(\zeta) = -\zeta+\pi\zeta^{2k+1} + o(\zeta^{2k+1}),$$
hence the origin is the merging of a fixed point with $k$ periodic orbits of period 2, as in Example 2. 
At the level of the real analytic vector field, we have the merging of a singular point with $k$ limit cycles. 

\subsection{Example 5:} 
\begin{figure}\begin{center} \includegraphics[width=4cm]{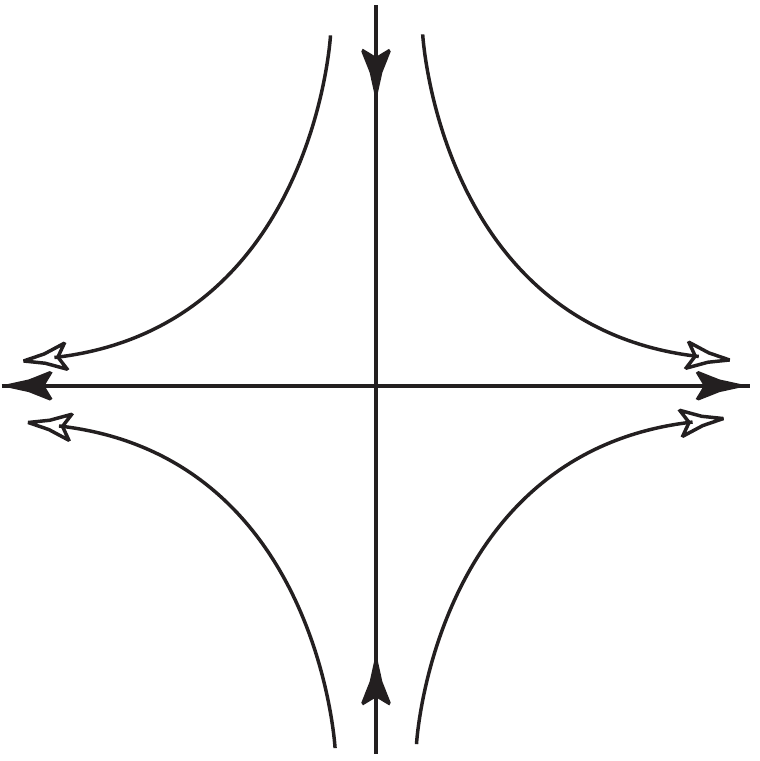}\caption{A resonant saddle.}\label{fig:saddle}\end{center}\end{figure}

A germ of resonant saddle of a planar vector field of
codimension $k$ with quotient of eigenvalues
$-\frac{p}{q}$ (see Figure~\ref{fig:saddle})
 and with orbital normal form
\begin{equation}\begin{cases}
\dot x = x,\\
\dot y = y\left(-\frac{p}{q} +x^{kp}y^{kq}+
ax^{2kp}y^{2kq}\right).\end{cases}\label{resonant_saddle}\end{equation}

As such the origin seems a simple hyperbolic point. To understand where multiplicity comes from, it helps revisiting Example 4. 
Taking $w=\overline{z}$, the system can be rewritten
$$\begin{cases} \dot z = iz + O(|z,w|^2),\\
\dot w = -iw + O(|z,w|^2).\end{cases}$$
Extending $z,w$ to $\C$, we have a resonant saddle of a complex singular foliation, corresponding orbitally to the case $p=q=1$ in \eqref {resonant_saddle}. And we have seen that we had the merging of $k$ limit cycles with the origin. When perturbing the vector field, the $k$ limit cycles lie on special leaves of the foliation that have non trivial homology. We will have the same behaviour for any germ of resonant saddle point. 

For such a point, the stable and unstable manifolds (which are also called the separatrices) are analytic and an analytic  change of coordinates can bring them to the coordinate axes. Similar to the case of the saddle-node (Example 3), we can define the holonomies of the two sepatrices of the vector field (see Figure~\ref{fig:holonomy}). Martinet and Ramis proved \cite{MR2}  that two systems with orbital formal normal form \eqref{resonant_saddle} are analytically equivalent if and only if the holonomies of their $x$-separatrices are conjugate and, by symmetry, if and only if the holonomies of their $y$-separatrices are conjugate. The holonomies $h_x$ and $h_y$ of the $x$- and $y$-separatrices respectively each have a resonant fixed point of codimension $k$ at the origin of the form of Example 2:
$$\begin{cases}h_x(y) = \exp\left (-\frac{2\pi i p}{q}\right)y + o(y),\\
h_y(x) =  \exp\left (-\frac{2\pi i q}{p}\right)x + o(x).
\end{cases}$$
Moreover,  Yoccoz and P\'erez-Marco \cite{PMY} showed that any $1$-dimensional germ of analytic diffeormorphism with a fixed point can be realized as a holonomy map of a germ of saddle point of an analytic  vector field in $\C^2$.

\subsection{Example 6:} 
\begin{figure} \begin{center}\includegraphics[width=6cm]{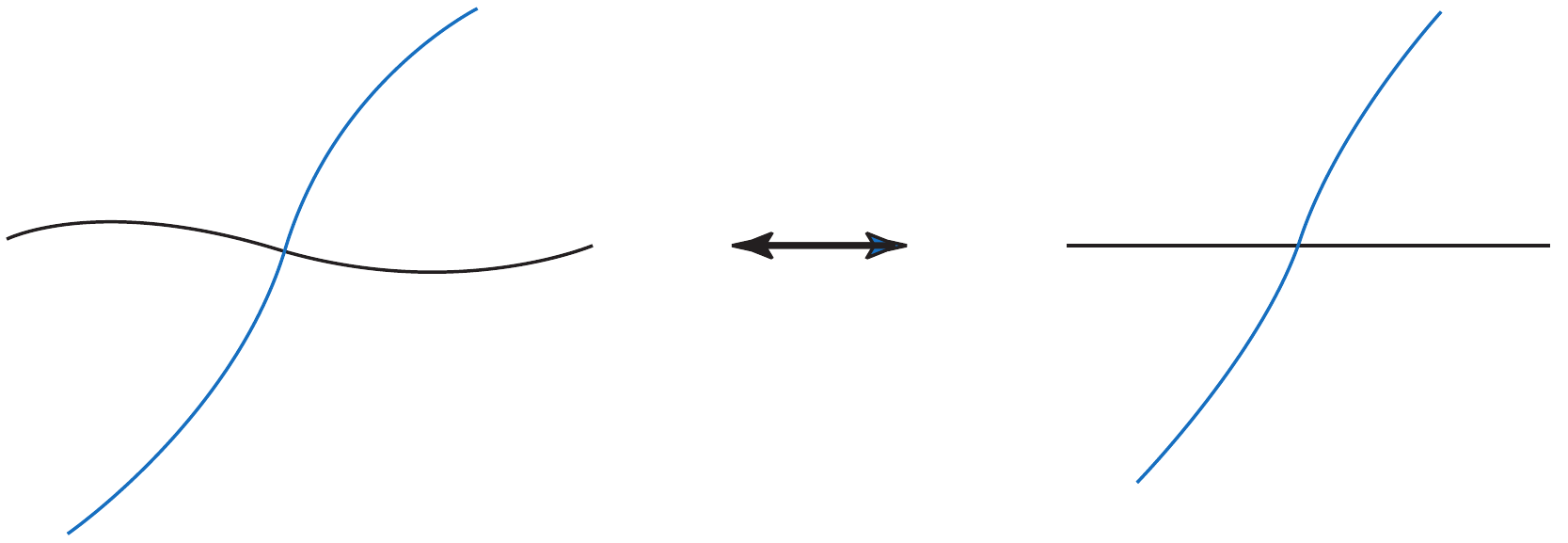}\caption{Conformal equivalence of curvilinear angles.}\label{fig:curvilinear}\end{center}\end{figure}

A germ of curvilinear angle  formed by two germs of real analytic curves (see Figure~\ref{fig:curvilinear}), where the angle is of the form
$2\pi\frac{p}{q}$, which we call a \emph{rational angle}.

In this case the equivalence problem can be stated: \emph{ When are two germs of curvilinear angles
conformally equivalent, i.e. mapped one onto the other by a germ of conformal diffeormorphism?}

 A  
Schwarz symmetry  $z\mapsto \Sigma_j(z)$ is associated to each curve of the angle. (The Schwarz symmetry $\Sigma$ of a germ $\gamma$ of real analytic curve is defined as follows: let $h$ be a germ of holomorphic diffeomorphism sending $\gamma$ to the real axis. Then $\Sigma=h^{-1}\circ \sigma\circ h$, where $\sigma(z) = \overline{z}$.)

Let $$f=\Sigma_2\circ \Sigma_1.$$

Then $f$ is a germ of analytic
diffeomorphism
$f(z) = e^{4\pi i\frac{p}{q}}z + o(z),$ 
satisfying \begin{equation}\Sigma_1\circ f= f^{-1}\circ
\Sigma_1.\label{cond:curvi}\end{equation} Let us call $f$ the \emph{associated diffeomorphism to the curvilinear angle}. Then $f$  is of the type studied in Examples 1 and 2, together with an additional symmetry. 
When considering the conformal equivalence of two curvilinear angles corresponding to Schwarz symmetries $\Sigma_1, \Sigma_2$ for the first angle, and $\Sigma_1', \Sigma_2'$ for the second angle, we can of course suppose that $\Sigma_1=\Sigma_1'=\sigma$. Then two germs of curvilinear angles are conformally equivalent if and only if their associated diffeomorphisms are conjugate under a change of coordinate $h$ commuting with $\sigma$:\
$$ f'= h\circ f\circ h^{-1}, \qquad h\circ \sigma=\sigma\circ h.$$ 

\subsection{Example 7:} A nonresonant irregular singular point of Poincar\'e rank $k$
of a linear differential system

$$x^{k+1}\frac{dy}{dx} = A(x) y =(D_0+O(x))y, \qquad y\in \C^n,$$
with normal form $$x^{k+1}\frac{dy}{dx} = (D_0+ D_1x+ \dots + D_k x^k)
y,$$
where  $D_0,D_1, \dots D_k$ are diagonal and the eigenvalues of $D_0$ are distinct (this is the nonresonance condition.)
The singular point at the origin has multiplicity $k+1$. This example is a bit different from the first six: there is no one dimensional map to explain the dynamics. It will be discussed in detail in Section~\ref{sec:Example7}.

\section{The example of the saddle-node}\label{sec:sn}

The formal normal form at a saddle-node is given by  
\begin{align}\begin{split} 
\dot x&=x^2,\\
\dot y&=y (1+Ax).\end{split} \label{nf_saddle-node}\end{align}

\bigskip

\noindent{\bf Questions.}

\begin{enumerate}
\item We cannot get rid of the term $Axy$ in the second equation. What is the meaning of  $A$?
\item Generically the change to normal form diverges: {\bf  Why?}\end{enumerate}

\bigskip
\noindent{\bf Spirit of the answer to the second question:} 

\begin{enumerate} 
\item We need to extend $(x,y)$ to a small polydisk in  $\mathbb C^2$.
\item  The saddle-node is a multiple singular point. Hence it is natural to unfold it. The orbital formal normal form of a generic unfolding is \begin{align}\begin{split} 
\dot x&=x^2-\eps,\\
\dot y&=y (1+A(\eps)x).\end{split}\label{sn_unfold} \end{align}
In the unfolding, there are rigid \emph{models} near each of the two singular points.
Generically these models mismatch till the merging of the singular points, yielding divergence at the limit.\end{enumerate}

\begin{example} {\bf One example of mismatch.}
Consider the unfolding of a saddle-node with normal form \eqref{nf_saddle-node}
It is known that generically a saddle-node has no analytic center manifold. A famous equation is 
\begin{align}\begin{split} 
\dot x&=x^2,\\
\dot y&=-x^2+y,\end{split}\label{no_anal_cm} \end{align}
where the formal series for the center manifold is the divergent series $$y=\sum_{n\geq1} (n-1)!\:x^{n+1},$$ which is the solution asymptotic to $0$ of the linear differential equation $x^2y'-y+x^2=0$. 
This divergent series is \emph{$1$-summable} and its sum yields a function defined over a sector of opening $3\pi$: $\arg x\in \left(-\frac{\pi}2,\frac{5\pi}2\right)$ (see Figure~\ref{fig3}(a)).

In the general case of a germ of saddle-node of the form\begin{align}\begin{split} 
\dot x&=x^2,\\
\dot y&=f(x)+y(1+Ax+O(x^2)) +O(y^2),\end{split}\label{no_anal_cm2} \end{align}
where $f(x) =O(x^2)$, the formal series for the center manifold would also be generically divergent and $1$-summable with a sum defined on a sector
$$V=\left\{x\: :\: |x|<r, \arg{x}\in \left(-\frac{\pi}2+\delta,\frac{5\pi}2-\delta\right)\right\}$$ with $\delta>0$ arbitrarily small and $r$ depending on $\delta$ (in the general case $r\to0$ when $\delta\to 0$). 

\begin{figure}[h]
\begin{center}
\subfigure[ $\eps=0$]
    {\includegraphics[angle=0,height=3.2cm]{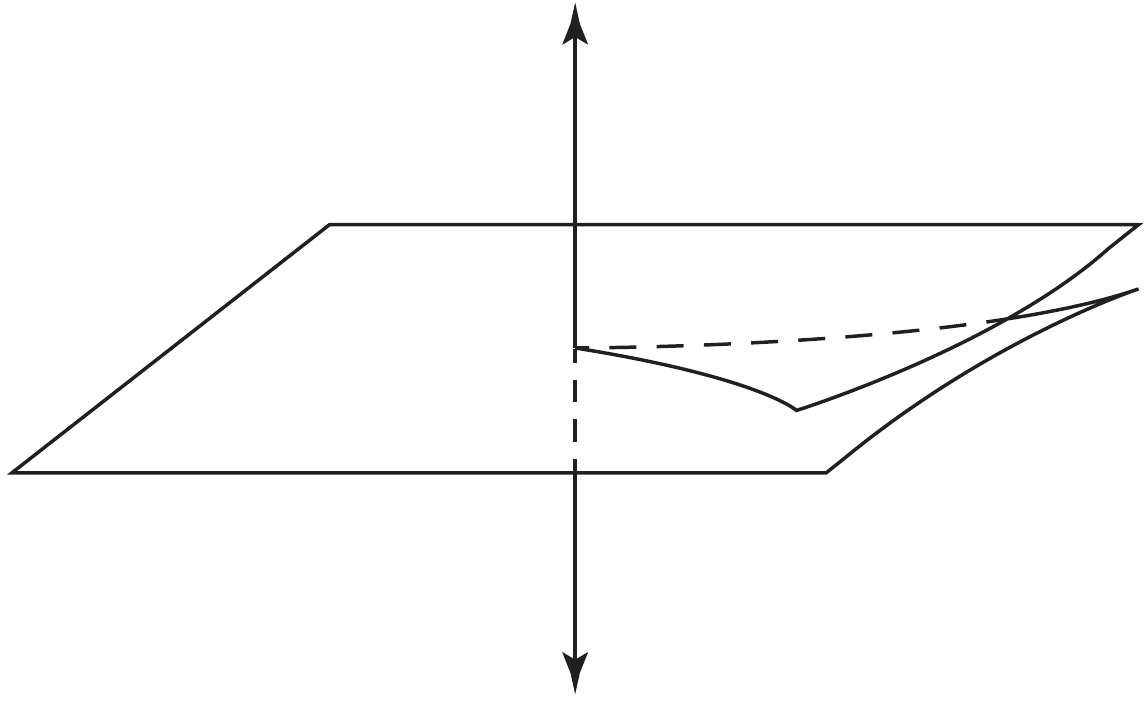}}
     \qquad
\subfigure[ $\eps\neq0$]
    {\includegraphics[angle=0,height=3.2cm]{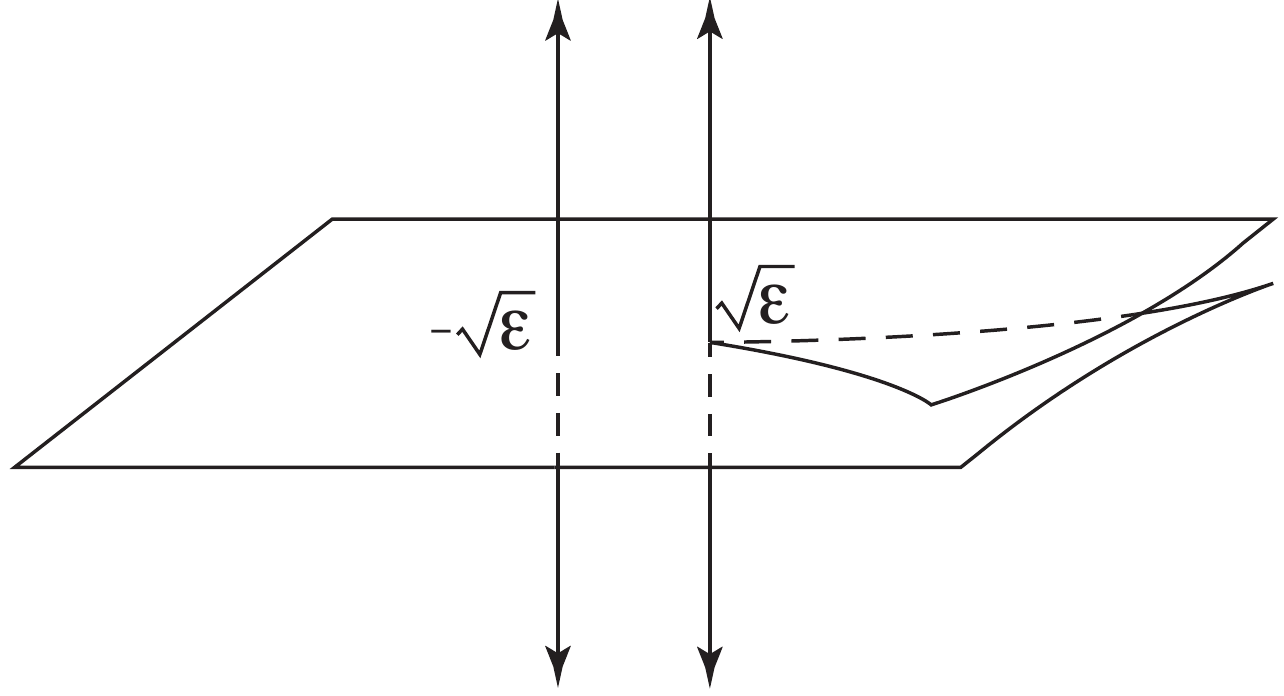}}
\caption{The generic form of a center manifold: we have drawn $x\in \C$ and $y\in \R$.}\label{fig3}
\end{center}
\end{figure}

This could seem mysterious. It is not, if we consider an unfolding
\begin{align}\begin{split} 
\dot x&=x^2 -\eps,\\
\dot y&=f_\eps(x)+y(1+A(\eps)x+O(x^2)) +O(y^2),\end{split}\label{no_anal_cm3} \end{align}
where $f_\eps(x) = O(x^2-\eps)$.
Indeed, for $\eps>0$ we now have two singular points, a saddle at $(-\sqrt{\eps},0)$ and a node at $(\sqrt{\eps},0)$. The saddle has a unique analytic stable manifold and the saddle-node has a unique center manifold on the side $x<0$.  
Hence it makes sense that the center manifold is the limit of the stable manifold of the saddle when $\eps\to 0$. 
Let us now look at the node. It has one small eigenvalue $\lambda_1= 2\sqrt{\eps}$, and a large eigenvalue $\lambda_2=1+O(\sqrt{\eps})$. When the node is nonresonant, i.e. $\lambda_2/\lambda_1\notin\N$, then the node is linearizable, by an analytic change of coordinates $(x,y)\mapsto (X,Y)$. Apart from the strong manifold $X=0$, the solutions in the neighborhood of the node are of the form $Y=CX^{\lambda_2/\lambda_1}$ (see Figure~\ref{fig:node}). As complex functions they are all ramified except the one for $C=0$, which we call the \emph{weak analytic manifold}. Hence, among all the trajectories with $\alpha$-limit at the node, the node has two exceptional ones which are analytic.
Suppose that the center manifold is not analytic. Then, for small $\eps>0$ we have no choice: the analytic stable manifold of the saddle must be ramified at the node (see Figure~\ref{fig3}(b)). 
This means that it does not coincide with the weak analytic manifold of the node. This is of course the generic situation. If this situation persists until the limit, then we are not surprised that we may have a ramified center manifold, and even that this situation should be the generic one. The divergence comes from the limit of this mismatch of the two local models at the saddle and at the node, each having a distinguished weak manifold. We have reached our:
\begin{figure}\begin{center}\includegraphics[width=4cm]{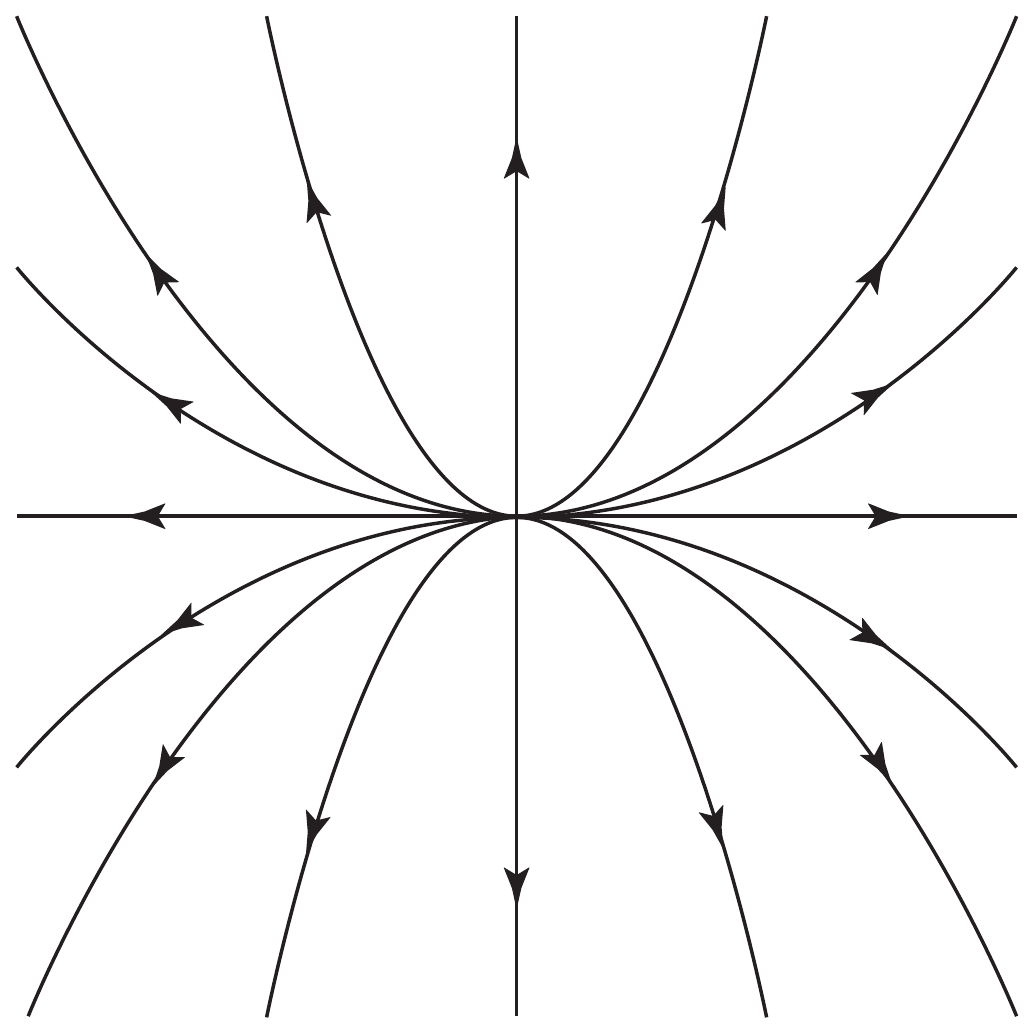}\caption{A node.}\label{fig:node}\end{center} \end{figure}

\medskip

\noindent{\bf Conclusion 1:} When we unfold a system
with no analytic center manifold, then the
analytic separatrices of the two singular points do not match.

\medskip Let us now consider the sequence of parameter values $\eps_n\to 0$ for which the node is resonant, i.e. $\left.\lambda_2/\lambda_1\right|_{\eps= \eps_n}= n\in\N$. 
In this case, the node is analytically conjugate  by a change of coordinates $(x,y)\mapsto (X,Y)$ to its normal form
\begin{align*}\begin{split}
\dot X &= \lambda_1 X,\\
\dot Y &= \lambda_2 Y +B_nX^n.\end{split}\end{align*} 
If $B_n=0$, then all solutions apart from the strong manifold have the form $Y=CX^n$, and none of them is ramified. Hence, this situation is forbidden when $\eps_n$ is sufficiently small, and the center manifold at $\eps=0$ is not analytic. We are then forced to have $B_n\neq0$, when $\eps_n$ is sufficiently small, i.e. $n$ sufficiently large. In that case all solutions apart from the strong manifold are ramified since of the form
$$Y=\frac{B}{\lambda_1}X^n\log X + CX^n.$$
That means that on this convergent sequence $\{\eps_n\}$ of parameter values the mismatch is carried by the singular point itself. We have reached our:

\medskip

\noindent{\bf Conclusion 2:} When we unfold a system
with no analytic center manifold then the node
is non linearizable as soon as resonant. This is the  \emph{parametric resurgence phenomenon}.

\medskip 

We have understood why divergence of the center manifold is the generic case and convergence is the exception.
This is just one example of the many \emph{mismatches} that occur within analytic dynamical systems and that follow very similar rules to Conclusions 1 and 2. 
\end{example}

\medskip

And what about the formal invariant and Question 1?
Let us look again at the formal orbital normal form \eqref{sn_unfold} of the unfolding.   
The ratios of eigenvalues at the singular points are
$$\mu_\pm= \pm\frac{2\sqrt{\eps}}{1\pm A\sqrt{\eps}}$$
Then $$\frac1{\mu_+} + \frac1{\mu_-}= A,$$
which means that $A$ is a measure of the mismatch of these ratios of eigenvalues until the limit. 
In the \emph{orbital} formal  normal form, we have two singular points, each with its ratio of eigenvalues, hence two ratios of eigenvalues. 
Therefore, to allow full generality, we need two independent parameters to control these, namely $\eps$ and $A$.
We get our third conclusion:

 \medskip

\noindent{\bf Conclusion 3:} In the formal normal form of the unfolding, the number of parameters is equal to the number of analytic invariants at the linear level at each simple singular point. 

\medskip Again this is a very general rule in all the examples we are considering. For planar vector fields, if we would consider full normal forms and not only orbital normal forms, the rule would still hold. For instance the formal normal form of the unfolding of a saddle-node is 
\begin{align*}\begin{split} 
\dot x&=(x^2 -\eps)(B_0+B_1x),\\
\dot y&=y(1+Ax)(B_0+B_1x)\end{split} \end{align*}
The four parameters $\eps, A, B_0$ and $B_1$ allow for the four eigenvalues (two at each singular point) to behave independently.

 \medskip
 
 If we now look again at the orbital formal normal form \eqref{sn_unfold} of the saddle-node we see that 
\begin{equation}\frac1{\sqrt{\eps}}= \frac1{\mu_+}-\frac1{\mu_-}.\label{canonical_parameter}\end{equation} Hence we reach our last conclusion: 
 
 \bigskip

\noindent{\bf Conclusion 4:} The parameter of the formal normal form is an analytic invariant. 

\bigskip Again, this is a general feature of all the examples we are considering except Example 7.

\section{Revisiting the seven examples and going further}\label{sec:common}

\subsection{\bf The common features} 

Here are common features of all seven examples, some of which we have already seen, and some we will explore further. 
\begin{enumerate}
\item In each case we have coallescence of $k+1$ \lq\lq special objects\rq\rq, which come with their local
model. The special objects could be singular points, periodic orbits, limit cycles, special leaves. 
 \item To understand why we have divergence, we unfold: 
when we have $k+1$ objects (codimension $k$) this  leads to $k$-summability in the limit.
\item In the unfolding,  generically the divergence can be seen as the limit of
the gluing of the $k+1$ local models which are rigid. Hence, {\bf divergence is the rule and convergence is
exceptional.}
\item Except in Example 7, the parameters in the formal normal form of the unfolding are canonical, since they are
analytic invariants (up to an action of the rotation group of order $k$.) 
 \item In all cases we have a finite
parameter family representing a formal normal form:
\lq\lq the model family\rq\rq. This is because all resonances are consequences of a unique one: this is the $1$-resonance case. 
\item The extra formal parameter(s) are present to match, in the unfolding, 
the need of independent multipliers in the diffeomorphism case, or the need of independent eigenvalues (resp. ratios of eigenvalues) for vector fields, depending whether the equivalence relation is conjugacy (resp. orbital equivalence).
\item Except in Example 7, the \lq\lq dynamics\rq\rq\ can be reduced
to that of a $1$-dimensional map.
\item In all cases we
observe a \emph{parametric resurgence
phenomenon},  i.e. the unfolded singular points have
pathologies on discrete sequences of parameter values 
$\{\eps_n\}$ converging to $\eps=0$.
\item In all cases  the study of the dynamics of a generic $k$-parameter unfolding is governed by the dynamics on $\CP^1$ of the $k$-parameter polynomial vector field $\dot z = P_\eps(z)$ or $\dot z = Q_\eps(z)$, where
\begin{equation} P_\eps(z)=z^{k+1} + \eps_{k-1} z^{k-1} + \dots + \eps_1+\eps_0,\label{def:P}\end{equation} and \begin{equation} Q_\eps(z)=z^{k+1} + \eps_k z^k+ \eps_{k-1} z^{k-1} + \dots + \eps_1.\label{def:Q}\end{equation}The vector field $\dot z = P_\eps(z)$ is the generic unfolding of $\dot z= z^{k+1}$, and $\dot z = Q_\eps(z)$ is the generic unfolding preserving the singular point at the origin.
\end{enumerate}

\subsection{The first six examples revisited}
Let us now revisit our first six examples and discuss briefly the unfoldings in each case. 

\subsubsection{Example 1 revisited.}\label{sect:ex1_revisited} The parabolic point of codimension $k$  is the coallescence of $k+1$  fixed points.
Given a generic $k$-parameter unfolding of a parabolic germ, there exists a change of parameter such that its normal form is the time-one map
of 
\begin{equation}\dot z = \frac{P_\eps(z)}{1+a(\eps)z^k},\label{parabol_unfold_nf}\end{equation}
with $P_\eps(z)$ given in \eqref{def:P}.
Fixing this new parameter $\eps$, there exists a change of coordinate depending analytically on $\eps$ and bringing the unfolding to the \emph{prepared form}    $$f_\eps(z) =
z+P_\eps(z)\left(1+R_\eps(z)+P_\eps(z)h(\eps,z)\right),$$ 
where 
\begin{itemize} 
\item $R_\eps(z)$ is a polynomial in $z$ of degree $\leq k$ depending analytically on $\eps$, 
\item $h$ is an analytic function of $(z,\eps)$,
\item if $z_1, \dots z_{k+1}$ are the zeroes of $P_\eps$ (i.e. the fixed points of $f_\eps$), then $f_\eps'(z_j)= \exp\left(P_\eps'(z_j)\right)$.\end{itemize}
The parameters are unique up to the action of the rotation group of order $k$ in \eqref{parabol_unfold_nf} (see Section~\ref{sec:can_par} below).
\bigskip

\begin{figure}\begin{center}
\subfigure[Left domain]{\includegraphics[width=3.5cm]{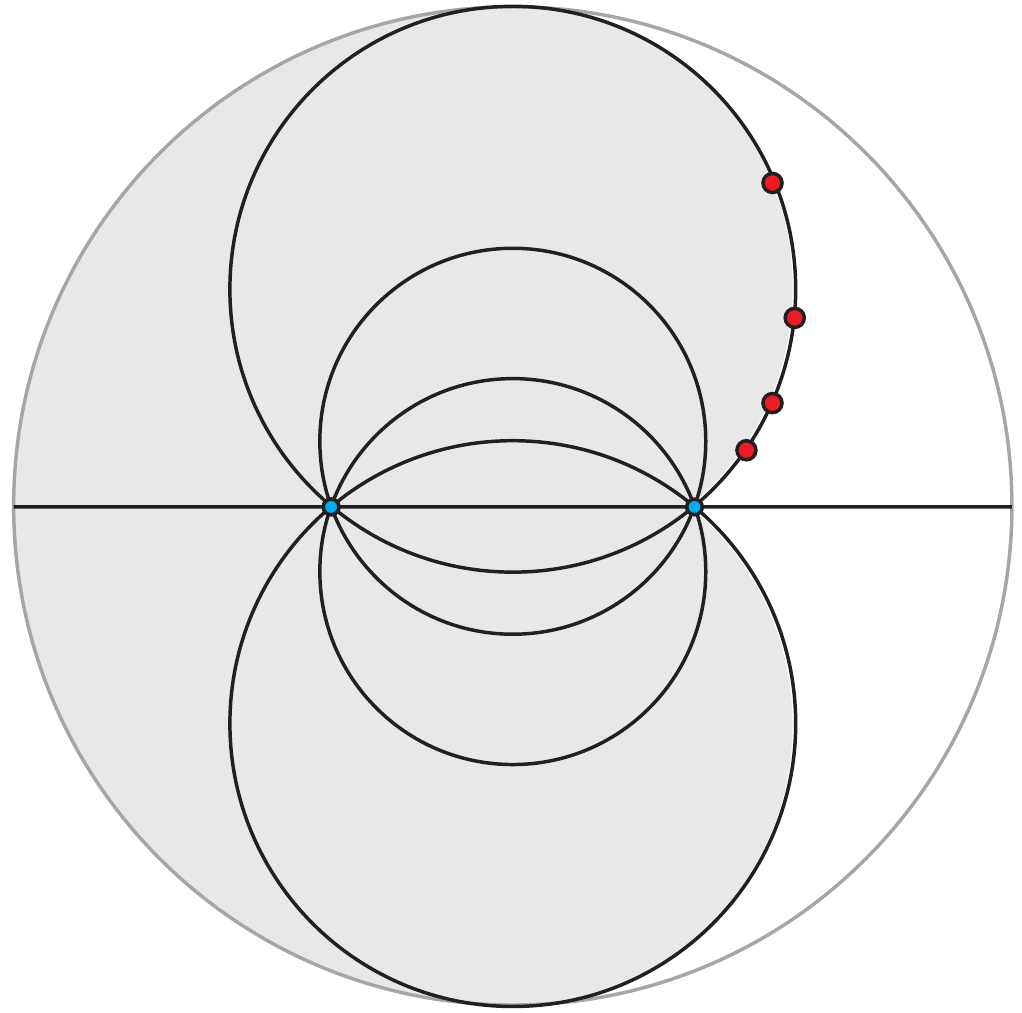}} \quad\subfigure[Right domain]{\includegraphics[width=3.5cm]{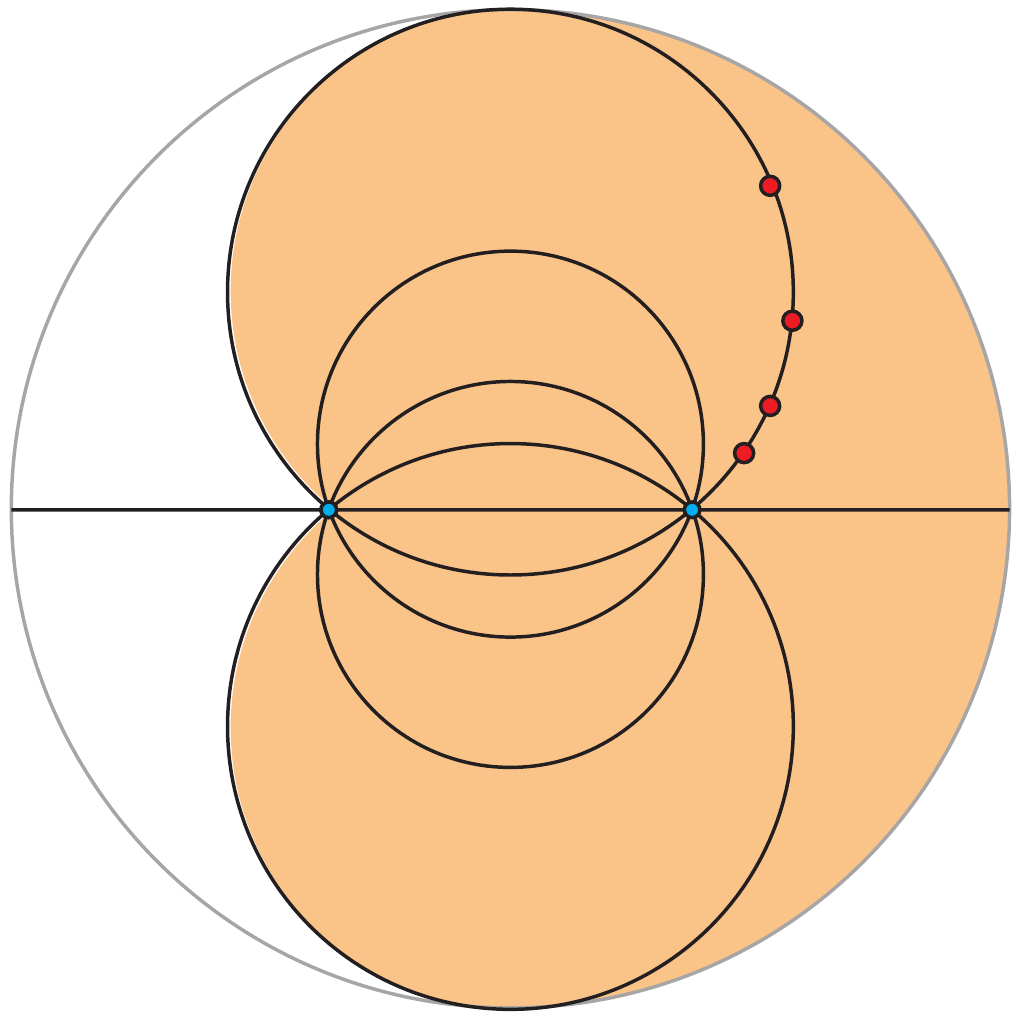}}\quad\subfigure[The intersection]{\includegraphics[width=3.5cm]{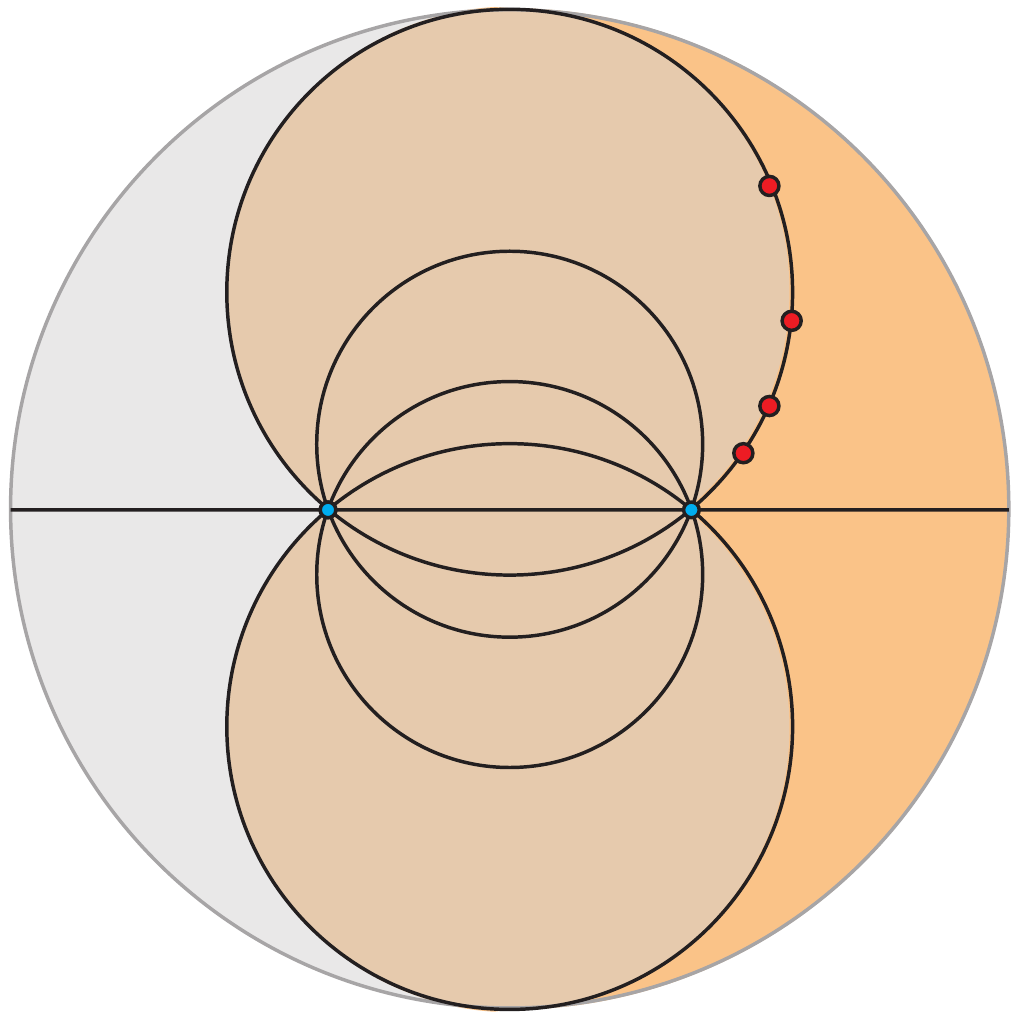}}\caption{For $\eps$ in the Poincar\'e domain, domains of normalizations near each fixed point intersect.}\label{fig:par_unf}\end{center}
\end{figure}

Let us now discuss the case $k=1$ in more detail. Two approaches were proposed in the literature. The first one was proposed, among others, by Arnold, and studied by Martinet \cite{Ma} and Glutsyuk \cite{Gl1}. It consists in unfolding in the \emph{Poincar\'e domain} in parameter space, namely a sector in parameter space where the fixed points are hyperbolic, one attracting, one repelling and, moreover, there are orbits going from the repelling fixed point to the attracting one. In that case, there exists an analytic change of coordinate to the normal form in the neighborhood of each fixed point, and the domains of normalization intersect (see Figure~\ref{fig:par_unf}). The normalizing changes of coordinates are unique up to a global symmetry of the formal normal form. A modulus of analytic classification comes from the comparison of these normalizing changes of coordinates. This approach was also considered for unfoldings of irregular singular points of linear differential equations by Ramis, Bolibruch and Glutsyuk among others (see for instance \cite{Ra} and \cite{Gl2}). For very long, no one in the field had any idea how to treat the other direction, since the fixed points could have multipliers given by an irrational rotation, a very difficult case, for which there is no modulus of analytic classification. 

\begin{figure}\begin{center}
\subfigure[Left domain]{\includegraphics[width=3.5cm]{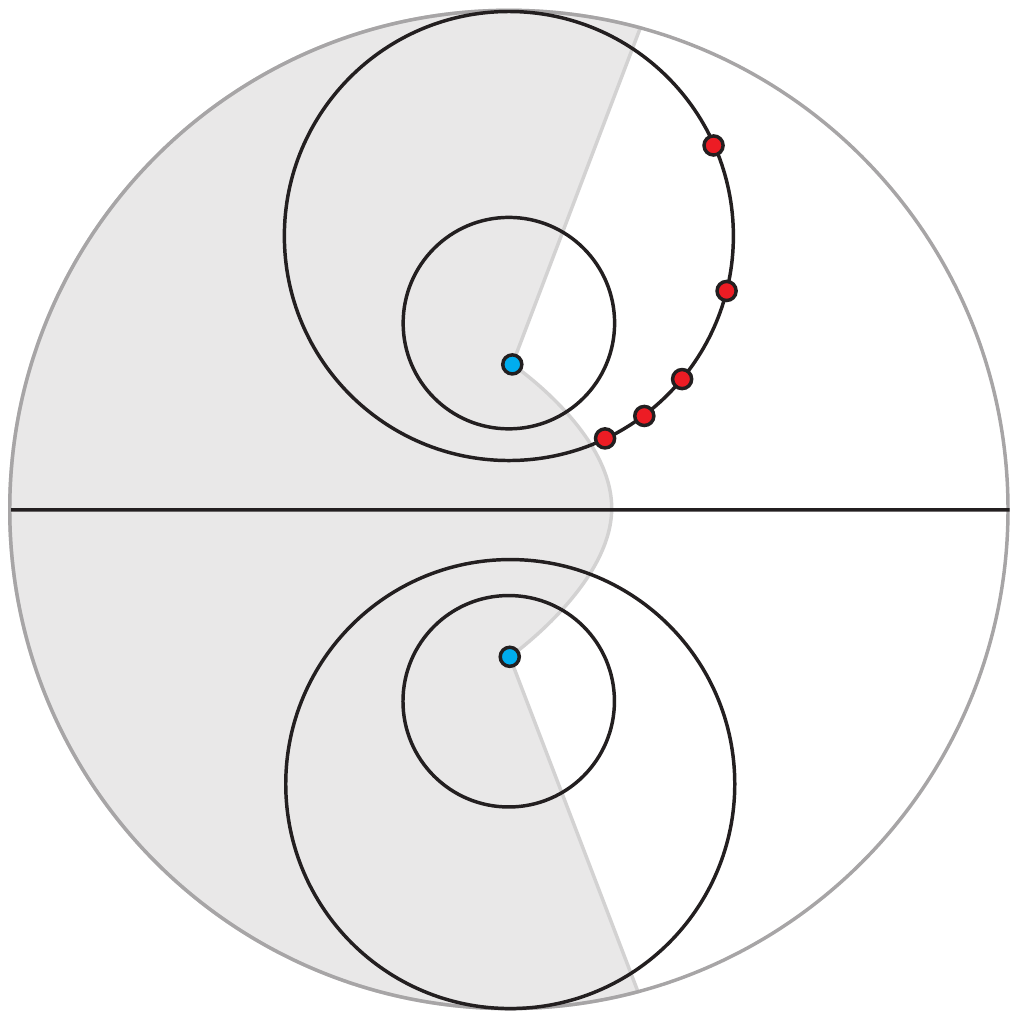}} \quad\subfigure[Right domain]{\includegraphics[width=3.5cm]{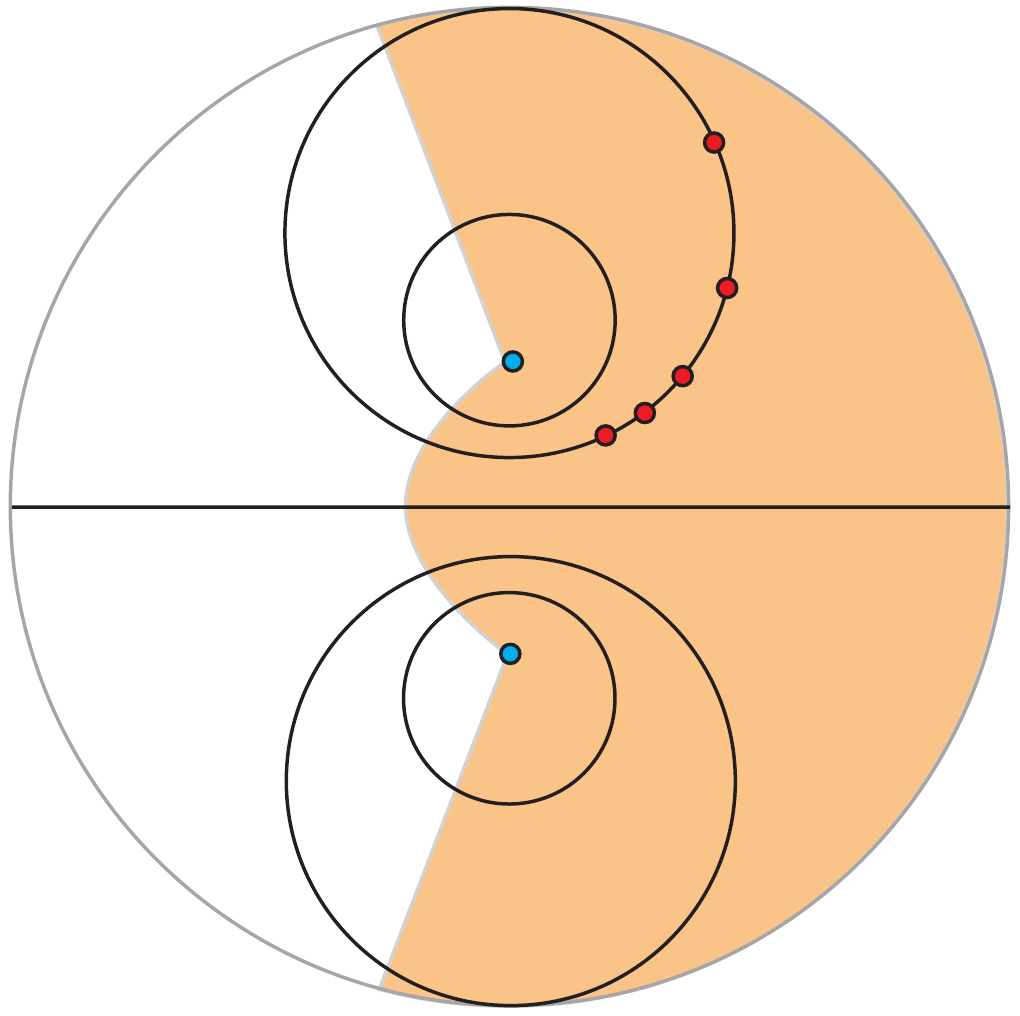}}\quad\subfigure[The intersection]{\includegraphics[width=3.5cm]{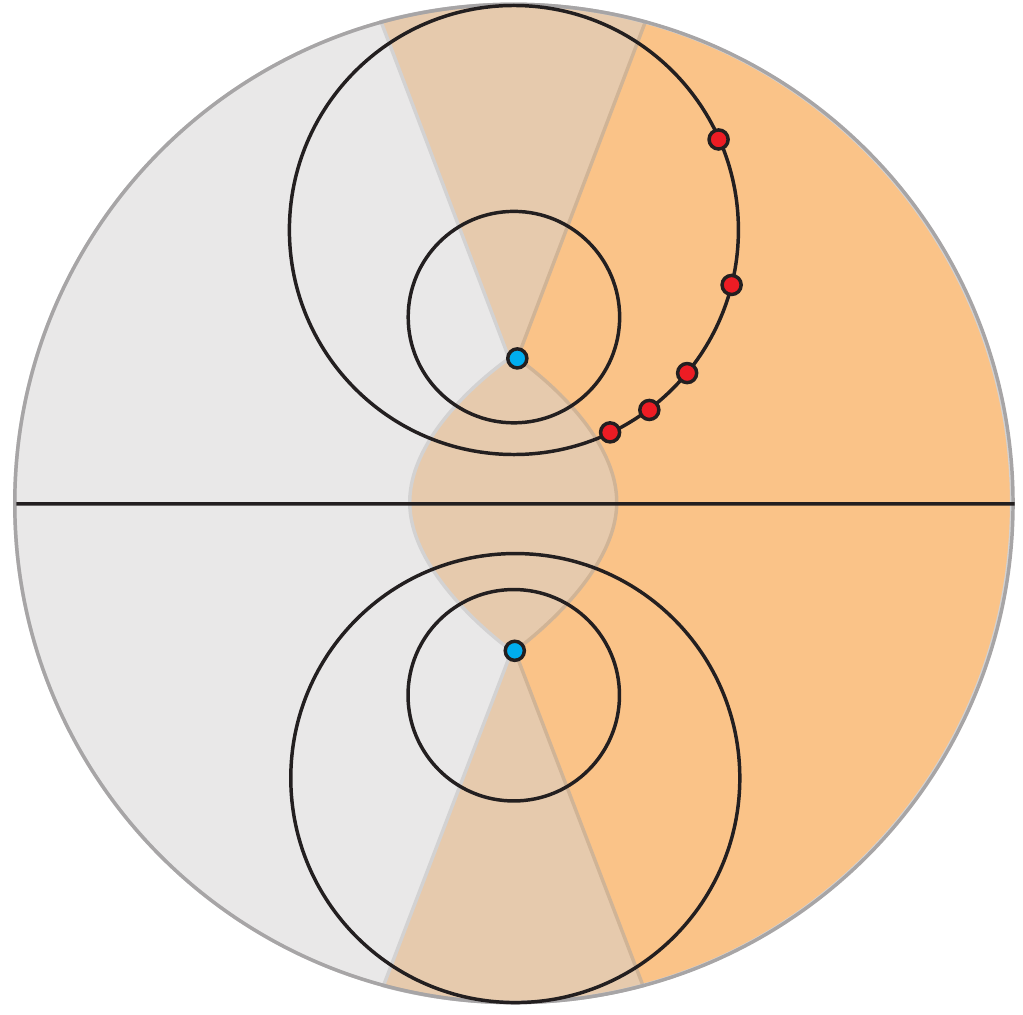}}\caption{For $\eps$ in the Siegel domain, domains of normalizations are adherent to the two fixed points.}\label{fig:domains_Lavaurs}\end{center}
\end{figure}

The second approach allows to treat the complementary sector, the \emph{Siegel domain}. It first appeared in the thesis of Lavaurs (\cite{La}), a student of Douady, when studying the parabolic implosion in Julia sets. 
As before, it consists in finding almost unique changes of coordinates to the normal form and to read the modulus by comparing these normalizations (see Figure~\ref{fig:domains_Lavaurs}).  But the difference is that the domains of these normalizing changes of coordinates have a sector shape in the neighborhood of each singular point. 

\bigskip
As in Section~\ref{sec:ex1}, we discuss the space of orbits in each case. 

\begin{figure} \begin{center} 
\subfigure[Poincar\'e domain]{\includegraphics[height=6cm]{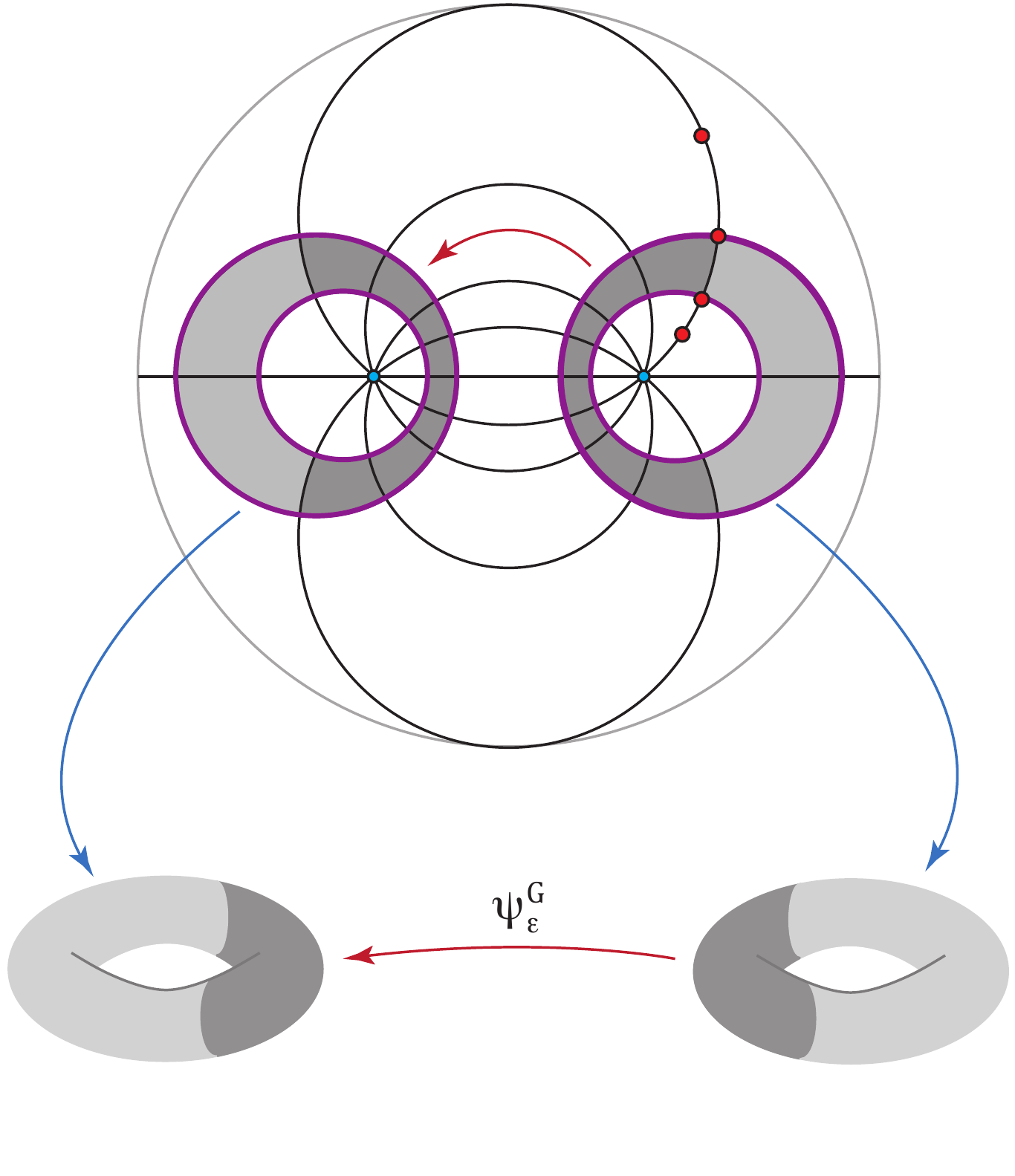}} \qquad\subfigure[Siegel domain]{\includegraphics[height=6cm]{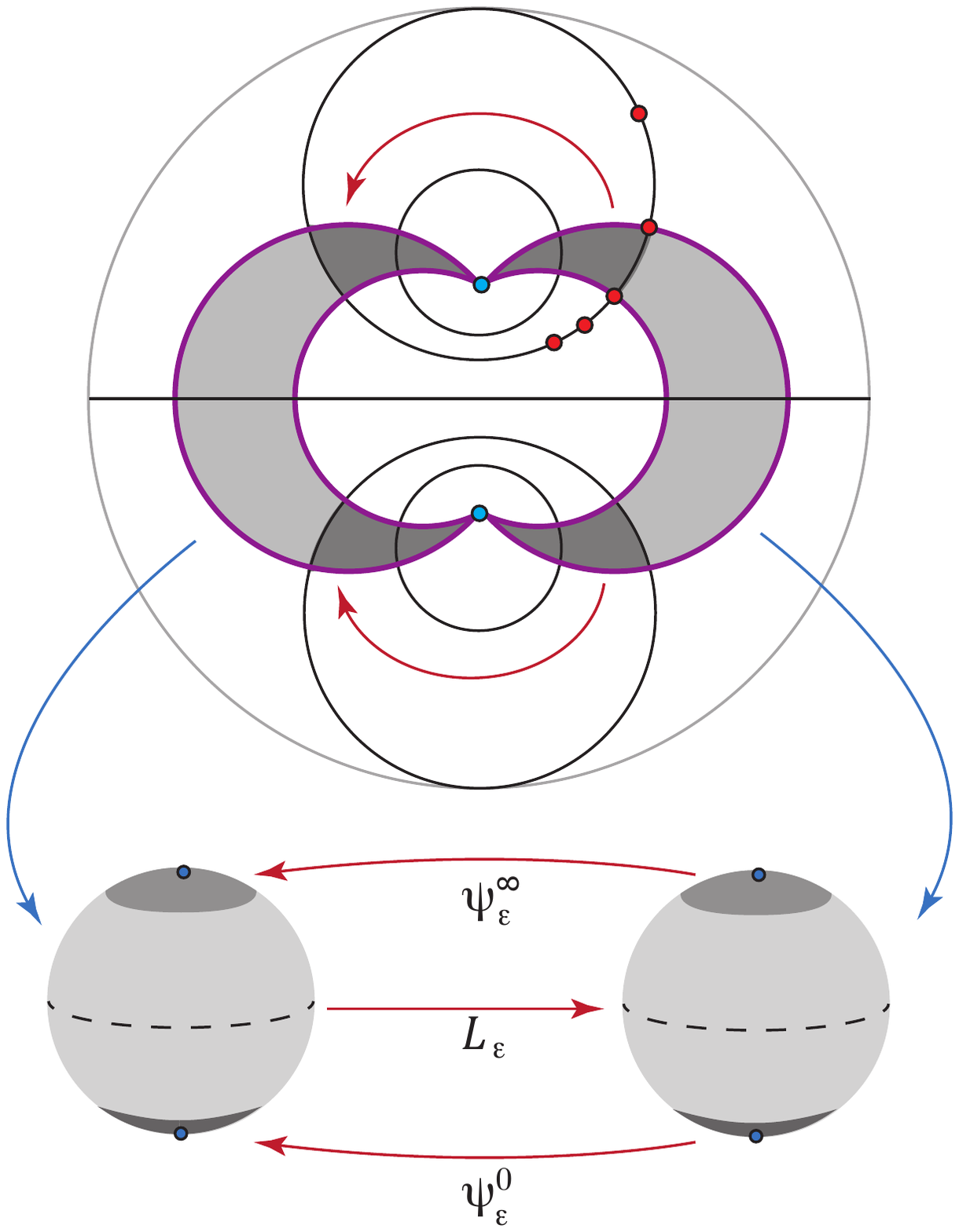}}\caption{The different ways of taking fundamental domains.}\label{parabolic_domain_unfolded}\end{center}\end{figure}

\bigskip\noindent{\bf Poincar\'e domain.} It this case, it is natural to take fundamental domains for each fixed point. These have the form of annuli: around each singular point we take a close curve $\ell$ sufficiently transversal to the flow lines of the vector field $\dot z = z^2-\eps$ so that  $f_\eps(\ell)$ is disjoint  from $\ell$ (see Figure~\ref{parabolic_domain_unfolded}(a)). Then $\ell$ and $f(\ell)$ bound an annular region. Using $f_\eps$ to identify $\ell$ and $f(\ell)$, each  fundamental domain is a torus. Hence, we have two tori corresponding to fundamental domains at each fixed point. The identification $\psi_\eps^G$ identifies two annuli, one on each torus. When $\eps\to 0$, each torus is pinched in the middle of its annulus and degenerates in a sphere with two points identified, thus cutting the annulus into two punctured disks. Hence, in the limit, the domain of  the transition map $\psi_\eps^G$ becomes disconnected.

\bigskip\noindent{\bf Siegel domain.} Here we take fundamental domains in the shape of crescents, with endpoints at the fixed points: the construction is the same as in the case $\eps=0$ in Section~\ref{sec:ex1} (see Figure~\ref{parabolic_domain_unfolded}(b)). Identifying the two sides of each crescent with $f_\eps$, yields fundamental domains conformally equivalent to spheres with two distinguished points $0$ and $\infty$ corresponding to the fixed points of $f_\eps$. As in the case $\eps=0$, this yields transition maps $\psi_\eps^0$ and $\psi_\eps^\infty$ defined respectively in the neighborhoods of $0$ and $\infty$ and identifying a point on one sphere to a point on the other sphere belonging to the same orbit. It is possible to choose coordinates on the spheres so that $\psi_\eps^0$ and $\psi_\eps^\infty$ depend analytically on $\eps\neq0$ with continuous limit at $\eps=0$.

 A new phenomenon occurs here: there is a global transition map $L_\eps$ between the two crescents, called the \emph{Lavaurs map}, which is linear when parametrizing the doubly punctered spheres with $\C^*$. 
Hence when $\eps\neq0$ belongs to the Siegel domain, one sphere is enough to describe the orbit space. In that case, the orbit space is given by one sphere quotiented by two  renormalized first return maps in the neighborhoods of $0$ and $\infty$, for instance $\tau_\eps^0= L_\eps\circ \psi_\eps^0$ and $\tau_\eps^\infty= L_\eps\circ \psi_\eps^\infty$ if we keep the right sphere. Note that the Lavaurs map has no limit when $\eps\to 0$, but it depends only on the formal normal form. Hence the first return maps also have no limit when $\eps\to 0$. But they are used to describe the parametric resurgence phenomenon. 

\bigskip\noindent{\bf The parametric resurgence phenomenon.} 
The phenomenon at the fixed point $-\sqrt{\eps}$ (resp. $\sqrt{\eps}$) is controlled by $\psi_\eps^0$ (resp. $\psi_\eps^\infty$). Let us discuss it at $-\sqrt{\eps}$. The Lavaurs map is of the form $w\mapsto K(\eps)w$, where $K(\eps)= \exp\left(\frac{C(\eps)}{\sqrt{\eps}}\right)$ for some nonzero $C(\eps)$ continuous, bounded, and bounded away from $0$.  
Hence, for each $p,q\in \N$ with $(p,q)=1$, there exists sequences of parameter values $\eps_n\to 0$ such that $K(\eps_n)(\psi_{\eps_n}^0)'(0)=\exp\left(2\pi i\frac{p}{q}\right)$. Since $\psi_\eps^0$ depends continously on $\eps$, then $K(0) = \lim_{n\to \infty}K(\eps_n)$ exists and $K(0)(\psi_0^0)'(0)=\exp\left(2\pi i\frac{p}{q}\right)$.
The \emph{parametric resurgence phenomenon} is simply the fact that if $K(0)\psi_0^0$ is non linearizable (because of some nonzero resonant  monomial in the formal normal form), then the renormalized first return map $\tau_{\eps_n}^0=K(\eps_n)\psi_{\eps_n}^0 = L_{\eps_n}\circ \psi_{\eps_n}^0$ is also nonlinearizable when $\eps_n$ is sufficiently small, i.e. $n$ is large enough. In the particular case $p=q=1$, this occurs as soon as $\psi_0^0$ is nonlinear. If the renormalized first return map is nonlinearizable at $0$, then this means that $f_\eps$ is nonlinearizable at $-\sqrt{\eps}$. Remember that the nonlinearity of $\psi_0^0$ is an obstruction to the convergence of the normalizing transformation. \emph{In the parametric resurgence phenomenon, the mismatch is carried by the fixed point itself.}

A similar phenomenon occurs at $\sqrt{\eps}$ when considering the renormalized first return map $L_\eps\circ \psi_\eps^\infty$ and adequate sequences of parameter values. 

\begin{remark}\label{par_resurg} \begin{enumerate}
\item Consider the parametric resurgence phenomenon under the hypothesis that $\tau_{\eps_n}^0=L_{\eps_n}\circ \psi_{\eps_n}^0$ is non linearizable of codimension $k$. Then $-\sqrt{\eps_n}$ is the coallescence of a fixed point with $k$ periodic orbits. Slightly perturbing slightly $\eps_n$ to $\eps'$ unfolds the situation and creates $k$ periodic orbits of $f_{\eps'}$ (multiplicity taken into account) around the fixed point. The larger $n$, the higher the period of the periodic orbits (all of the same period). 
\item 
Note that while $\tau_\eps^0$ has no limit when $\eps\to 0$, $\lim_{n\to \infty} \tau_{\eps_n}^0$ does exist.\end{enumerate} \end{remark}

\begin{figure} \begin{center} 
\includegraphics[height=9cm]{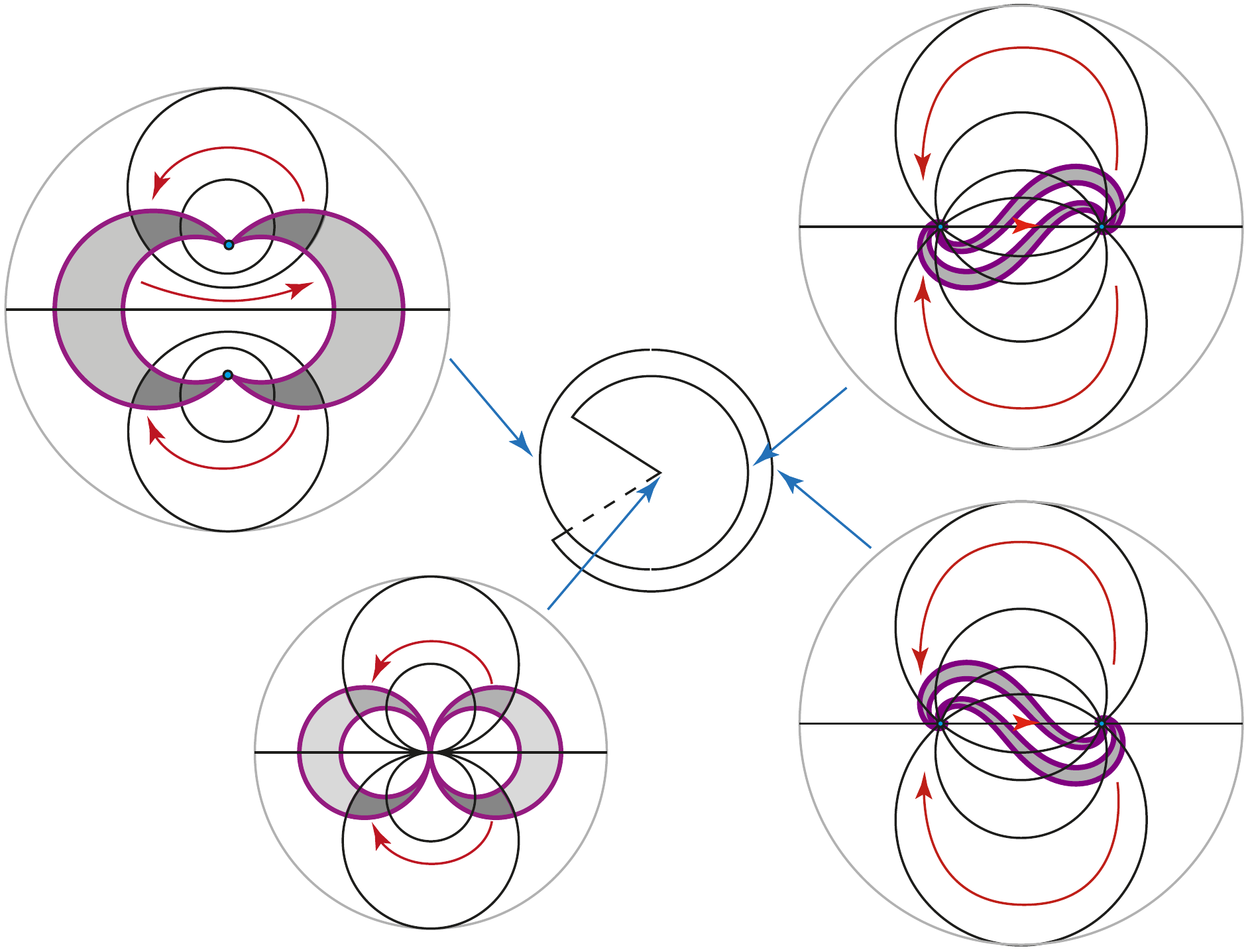} \caption{Lavaurs type fundamental domains for all parameter values.}\label{Bif_parabolic}\end{center}\end{figure}

\bigskip\noindent{\bf Covering the whole parameter space with the approach of the Siegel domain \cite{MRR}.} So far we have used two charts, the Poincar\'e domain and the Siegel domain, to cover the whole  parameter space. But this is not necessary. It is possible to cover the full parameter space by using fundamental domains in the shape of crescents from one fixed point to the other (see Figure~\ref{Bif_parabolic}). However, there is a double price to pay: 
\begin{itemize} \item
The crescents must spiral when the parameter is in  the Poincar\'e domain. 
\item The construction of continuous  fundamental domains is ramified in the parameter. \end{itemize}

\subsubsection{Example 2 revisited.} The germ of resonant fixed point of codimension $k$ corresponds to the coallescence of a fixed point with $k$ periodic orbits of
period $q$.
 Modulo a change of coordinate and parameter, a generic $k$-parameter  unfolding of $f$ is of the form
 $$f(z) = \exp\left(\frac{2\pi i p}{q} \right)z\left(1+\frac{1}{kq}Q_\eps(z^q)(1+O(z,\eps))\right),$$ 
 where 
 \begin{equation}Q_\eps(z)= z^k+ \eps_{k-1}z^{k-1} +\dots + \eps_1z+\eps_0.\label{def:Q}\end{equation}
And the normal form is the time-$1/q$ map of 
\begin{equation}\dot z = \frac{zQ_\eps(z^q)}{1+az^{kq}}\label{vf_Q}\end{equation} composed with the rotation of angle $\frac{2\pi p}{q}$. 
The parameters are unique up to the action induced by rotations of order $k$ in \eqref{vf_Q}.

This corresponds to a special slice of dimension $k$ in the parameter space of Example 1 revisited (Section~\ref{sect:ex1_revisited}) for the codimension $K=kq$.

\subsubsection{Example 3 revisited.} 
An orbital  formal normal form of a generic $k$-parameter unfolding of a saddle node of codimension $k$ is given by
\begin{align}\begin{split} 
\dot x&=P_\eps(x),\\
\dot y&=y \left(1+A(\eps)x^k\right),\end{split}\label{sn_unfold_k} \end{align}
for $P_\eps$ in \eqref{def:P}. Again we have $k+1$ parameters to control the $k+1$ quotients of eigenvalues and the parameters are unique up to the action of the rotations of order $k$.
The result of Martinet-Ramis can be generalized to the unfoldings: two unfoldings of a codimension $k$ saddle-node are orbitally equivalent if and only if the unfolded holonomies of the strong separatrices are conjugate. 

\subsubsection{Example 4 revisited.} 
\begin{figure}\begin{center}
\subfigure[$\eps<0$]
{\includegraphics[width=2.9cm]{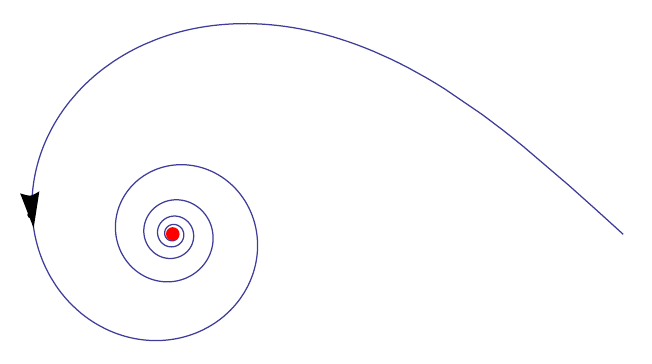}}\quad
\subfigure[$\eps=0$]
{\includegraphics[width=2.9cm]{Hopf_0}}\quad
\subfigure[$\eps>0$]
{\includegraphics[width=2.9cm]{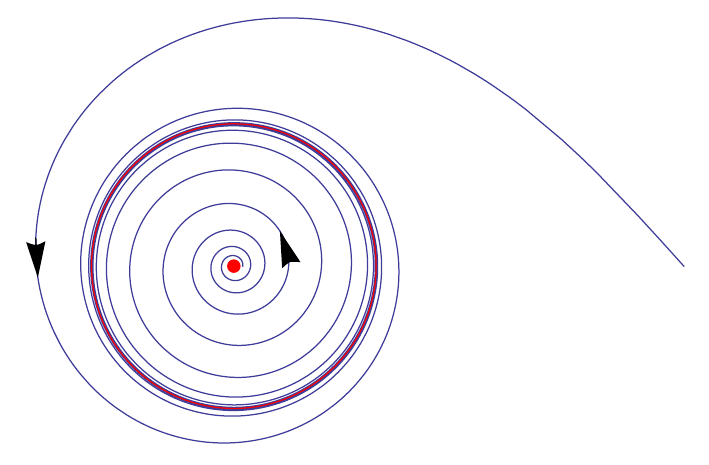}}\caption{The Hopf Bifurcation of order 1.}\label{fig:Hopf}
\end{center} \end{figure}

The generic $k$-parameter unfolding of a weak focus of order $k$ is the generic Hopf bifurcation of order $k$, which corresponds to the coallescence of a focus and $k$ limit cycles (see Figure~\ref{fig:Hopf} in the case $k=1$). It has orbital normal form
$$\dot z=
z\left(i+Q_\eps(z\overline{z})(1+a|z|^{2k})\right),$$
where $Q_\eps$ is given in \eqref{def:Q}.

Taking $w=
\ov{z}$, the orbital normal form can be rewritten
 \begin{align*} \begin{split}\dot z&=
z(i+Q_\eps(zw)(1+a(zw)^k)),\\
\dot w&=
w(-i+\ov{Q}_\eps(zw)(1+\ov{a}(zw)^k)),\end{split}\end{align*}
If we now take $(z,w)\in (\C^2,0)$, then this equivalent to the orbital unfolding of a resonant complex saddle. If $X_1, \dots, X_k$ are the zeroes of $Q(X)$, then each complex curve  $zw=X_j$ is a \emph{special  leaf} of the foliation, which has non trivial homology.

In \cite{AS} it is shown that two $1$-parameter generic real analytic families unfolding a weak focus of codimension $1$ are orbitally equivalent if and only if the families of their Poincar\'e return maps are conjugate. And, most probably,  the result is also true in higher codimension. 

\subsubsection{Example 5 revisited.} This is the case of the planar resonant saddle with orbital formal normal form \eqref{resonant_saddle}. 
An orbital  formal normal form of a generic $k$-parameter unfolding is 
$$\begin{cases}
\dot x = x,\\
\dot y = y\left(-\frac{p}{q}\left(1+Q_\eps(x^py^q)+ax^{2kp}y^{2kq}\right)\right),\end{cases}$$
where $Q_\eps$ is given in \eqref{def:Q}.
The result of Martinet-Ramis has been generalized in codimension $1$ (\cite{Ro3}) and is most probly true in higher codimension: 
two unfoldings of a codimension $k$ resonant saddles are orbitally equivalent if and only if the holonomies of a pair of corresponding  separatrices are conjugate. 
Generically, a $k$ parameter unfolding has $k$ special leaves with nontrivial holonomy which merge with the separatrices when $\eps= (\eps_0, \dots, \eps_{k-1})=0$. Hence, the $k+1$ \lq\lq special  objects\rq\rq, which merge together at $\eps=0$, are the $k$ special leaves together with  the singular point. 

\begin{remark}\label{rem:Yoccoz} This merging of the special leaves with the singular point is what has been called \emph{materialization of the Poincar\'e resonances} by Ilyashenko and Pjartli \cite{IP}. This phenomenon is very important. Indeed, consider a saddle point with an irrational ratio $-\alpha$ of eigenvalues. Then, in any neighborhood of $\alpha$, there are infinitely many rational numbers $p/q$. Generically, a special leaf (or special leaves) appears when the ratio of eigenvalues is perturbed from a rational value. Hence, it has been conjectured by Arnold and others that this accumulation of special leaves in the neighborhood of the singular point could be the obstruction to the linearizability of the saddle point when $\alpha$ is a Liouvillian irrational number (i.e. an irrational number well approximated by the rationals). 
On the other hand, if $\alpha$ is a Diophantian irrational number (i.e. badly approximated by the rationals), then the special leaves escape from a neighborhood of the origin sufficiently fast so that the point be orbitally linearizable. In view of \cite{PMY}, it suffices to consider fixed points of germs of $1$-dimensional diffeomorphisms, in which case the special leaves correspond to periodic points of the holonomy maps of separatrices. Yoccoz \cite{Y} proved that the mechanism of accumulation of periodic points in the neighborhood of a fixed point with multiplier $e^{i\pi \alpha}$  when $\alpha$ is irrational Liouvillian does indeed exist,  thus preventing linearizability, and that it is present in particular in the quadratic map $f(z) = e^{2\pi i \alpha}z + z^2$. Later, P\'erez-Marco (\cite{PM1} and \cite{PM2}) showed  that there exist other mechanisms preventing orbital linearizability of the fixed point: a fixed point with no periodic points in a neighborhood can be non linearizable.  \end{remark}

\subsubsection{Example 6 revisited.} 
\begin{figure}\begin{center}\subfigure[rational angle]{\includegraphics[width=5cm]{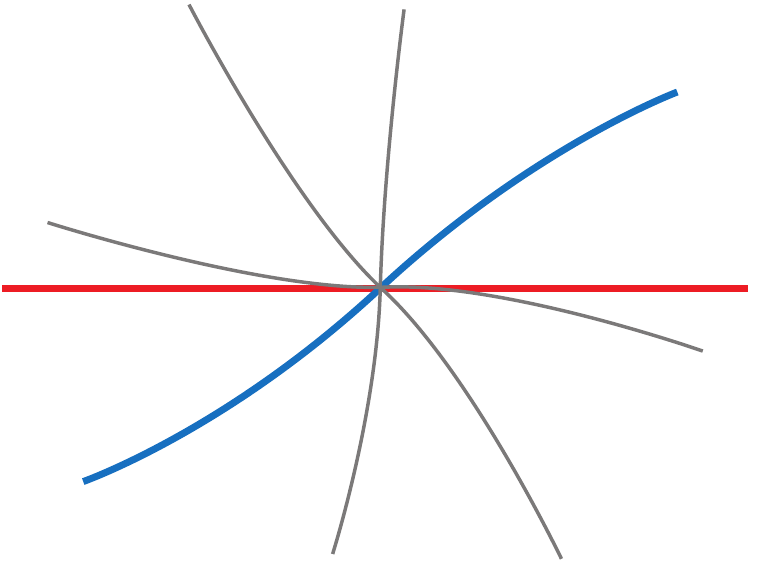}}\qquad\subfigure[deformation]{\includegraphics[width=5cm]{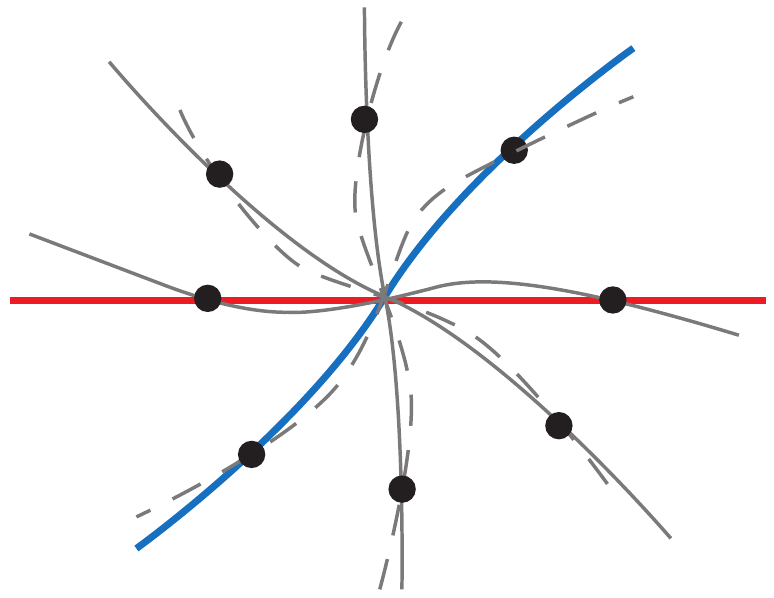}}\caption{Copies of the angle by symmetry with respect to the sides (a) until the two curves are tangent; (b) all copies pass through the periodic points after deformation. The dotted lines correspond to the 4-th, 5-th and 6-th copies of the angle.}\label{angle_copies}\end{center}\end{figure}

Let us look at a rational curvilinear angle of size $\frac{2\pi p}{q}$. Associated to each curve $\gamma_j$ of the angle is a 
Schwarz symmetry  $z\mapsto \Sigma_j(z)$  and we have introduced the analytic germ of diffeormorphism $f=\Sigma_2\circ \Sigma_1$ satisfying $f(z) = e^{4\pi i\frac{p}{q}}z + o(z).$ It could seem surprising that $f$ is generically not analytically linearizable. It is not. Indeed, we can make symmetric copies of the angle: for instance, the curve $\gamma_3$, which is the symmetric image of $\gamma_1$ under $\Sigma_2$ corresponds to the Schwarz symmetry $\Sigma_3$; and the diffeomorphism corresponding to the angle curvilinear formed by $\gamma_1$ and $\gamma_3$ is $f^{\circ 2}$. We iterate the construction by taking the symmetric of $\gamma_2$ with respect to $\gamma_3$, etc. until we get $\gamma_{n+1}$ which is tangent to $\gamma_1$ (see Figure~\ref{angle_copies}(a)). 
Of course, generically two tangent curves do not coallesce but have a multiplicity of intersection at the contact point. This multiplicity of intersection is an obstruction to the linearizability of $f$. Note that the diffeomorphism associated to the angle between $\gamma_1$ and $\gamma_{n+1}$ is $f^{\circ n}$.

\bigskip 
Let us now slightly perturb the angles between $\gamma_1$ and $\gamma_2$. Then we unfold the multiple fixed point of $f^{\circ n}$. But if the angle $\frac{2\pi p}{q}$ is not a multiple of $\pi$, then the perturbed $\gamma_1$ and $\gamma_2$ will still have a unique intersection point. Hence, the multiple fixed point of $f^{\circ n}$ corresponds to the merging of a fixed point of f with periodic points (see Figure~\ref{angle_copies}(b)). And, if we continue to make copies of the angle after $\gamma_{n+1}$, all the new curves  pass through the periodic points drawn on the figure. 

\begin{figure}\begin{center}
\subfigure[$\eps=0$]{\includegraphics[width=3.5cm]{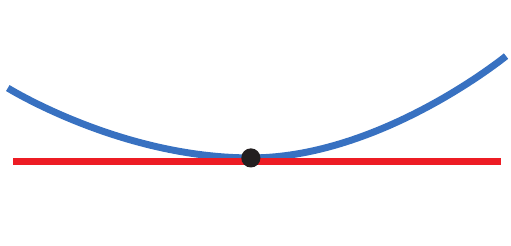}}\quad\subfigure[$\eps<0$]{\includegraphics[width=3.5cm]{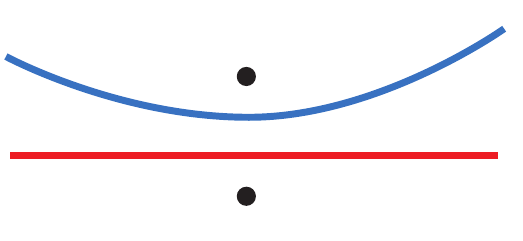}}\quad \subfigure[$\eps>0$]{\includegraphics[width=3.5cm]{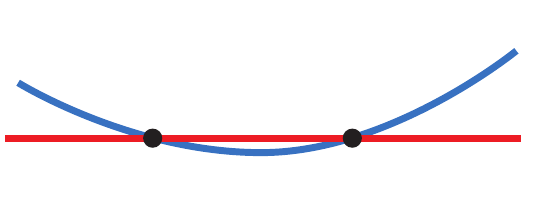}}\caption{The horn and its deformations in the two directions}\label{horn}\end{center}\end{figure}
\bigskip A particular case of curvilinear angle is that of the horn: two tangent curves (see Figure~\ref{horn}). In that case, in the normal form it is natural to take $\Sigma_1$ as the reflection with respect to the real axis. Then, because of \eqref{cond:curvi}, a formal normal form of the associated diffeomorphism  is the time-one of the vector field
$i\frac{z^2-\eps}{1+a(\eps)z}\,\frac{d}{dz}$. In the case $\eps<0$, the associated diffeomorphism $f_\eps$ has two fixed points outside the curve (Figure~\ref{horn}(b)). These points control the size of the neighborhood on which the unfolded angle can be analytically brought to the normal form: the neighborhood should not contain the points. 
In the case $\eps>0$ the two curves intersect at $\pm\sqrt{\eps}$ with (oriented) angles $\theta_\pm$ (Figure~\ref{horn}(c)). 
Then 
\begin{equation}\frac1{\theta_+} +\frac1{\theta_-}= a(\eps).\label{shift_angles}\end{equation}
Hence, $a$ introduces a difference between $|\theta_+|$ and $|\theta_-|$, which persists until the limit if $a(0)\neq0$. 

\begin{figure}\begin{center}\subfigure[$\eps=0$]{\includegraphics[height=5cm]{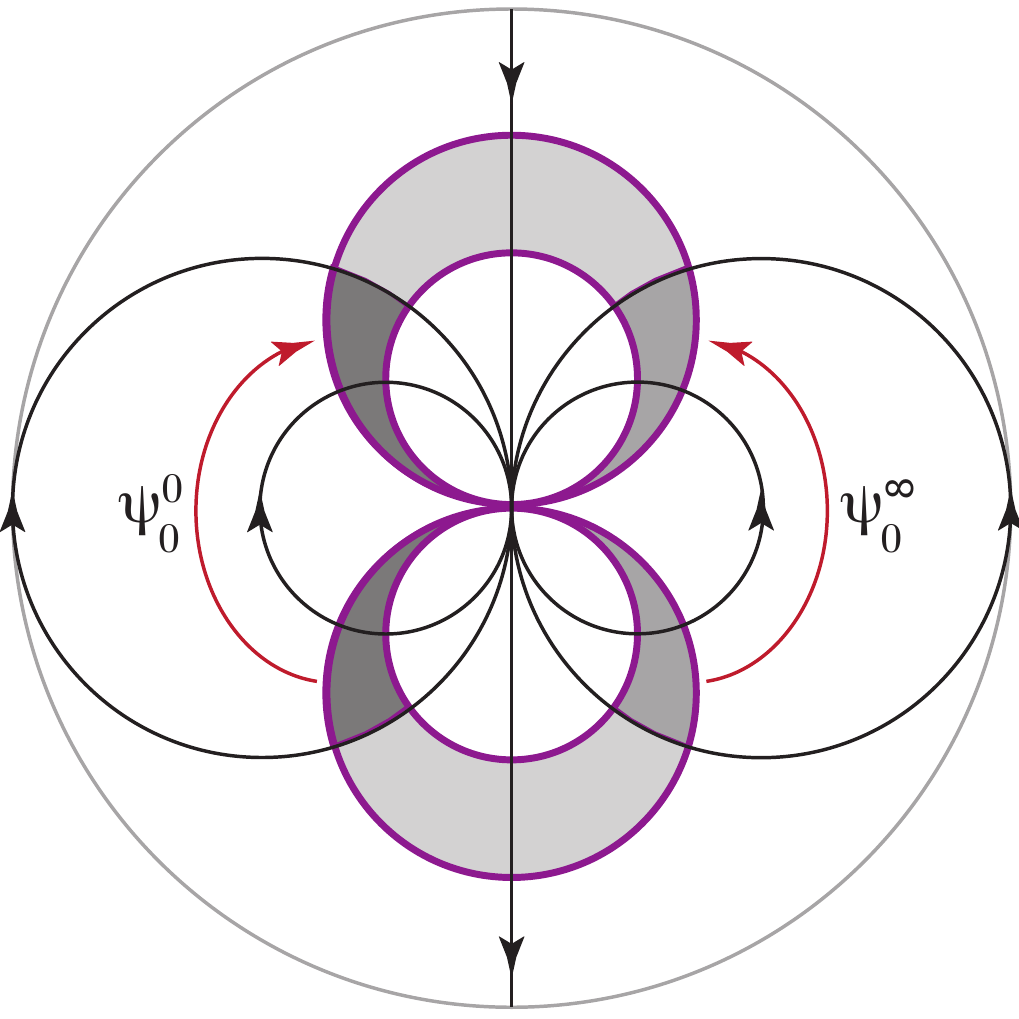}}\qquad\subfigure[$\eps>0$]{\includegraphics[height=5cm]{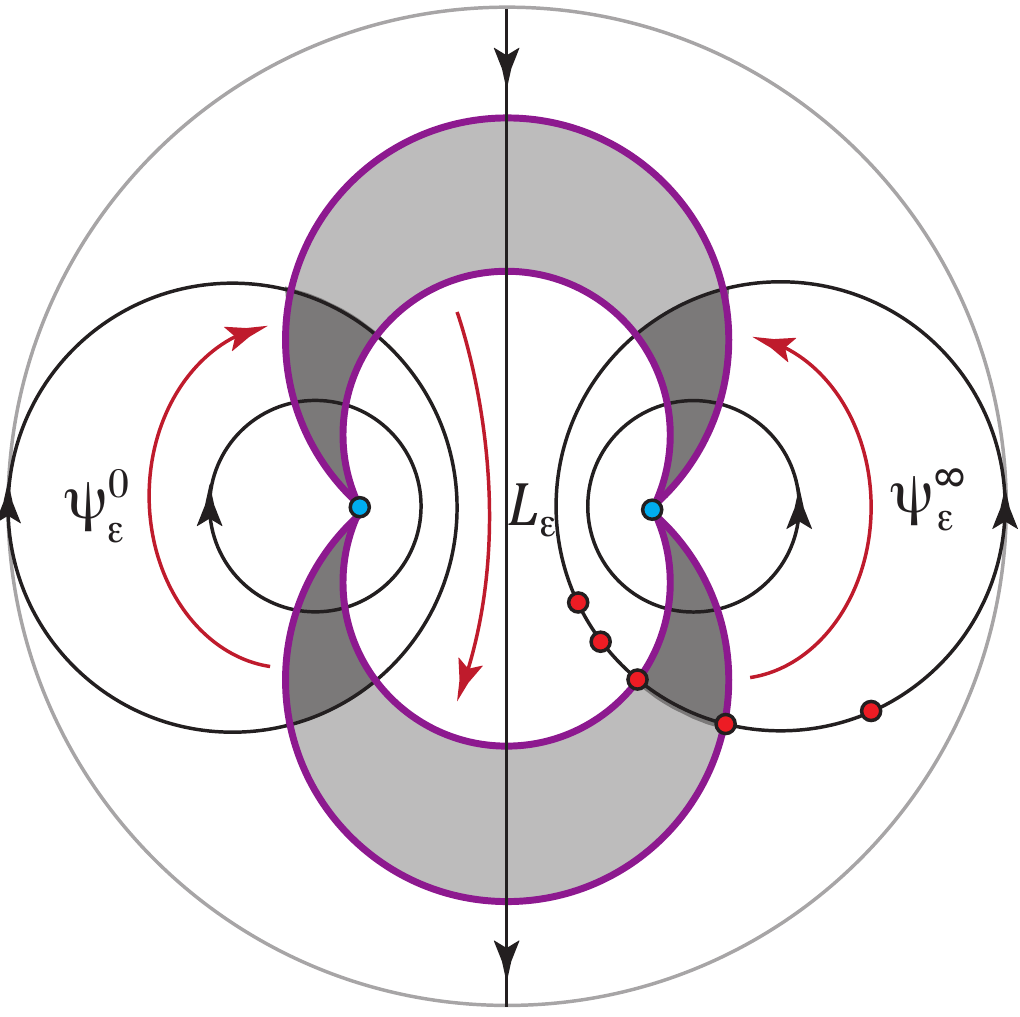}}\caption{The orbit space.}\label{orbit_space}\end{center}\end{figure}

\bigskip The space of orbits can of course be described as in Section~\ref{sec:ex1}: see Figure~\ref{orbit_space}. Let us now concentrate on the case $\eps>0$. Locally, at each intersection point we have a curvilinear angle and, depending on the value of $\eps$, this curvilinear angle may be rational, which means that we can apply to it the theory we just described. 
As explained in Section~\ref{sect:ex1_revisited}, the dynamics near $-\sqrt{\eps}$ (resp. $\sqrt{\eps}$) is described through the renormalized first return map, which has the form $\tau_\eps^0 = L_\eps\circ\psi_\eps^0$ (resp. $\tau_\eps^\infty= L_\eps\circ\psi_\eps^\infty$), and the phenomenon of parametric resurgence also occurs here.  

Let us discuss its meaning near $-\sqrt{\eps}$. Let $\{\eps_n\}$ be a sequence of parameter values for which $(\tau_{\eps_n}^0)'(0)= \exp\left(\frac{2\pi i p}{q}\right)$, which means that the curvilinear angle at $-\sqrt{\eps}$ is rational. Taking symmetric copies of this angle, we get a horn. The maximum multiplicity of intersection of the two curves forming this horn can be read from the codimension of the fixed point of $\tau^0 = \lim_{n\to \infty} \tau_{\eps_n}$ (see Remark~\ref{par_resurg}(2)). In the same way, the order of magnitude of the formal invariant (r\'esidu it\'eratif) of $\tau_{\eps_n}$ is given by that of $\tau^0$. And we have given an interpretation of the formal invariant as a kind of shift in the angles at two intersection points of two perturbed curvilinear curves (see \eqref{shift_angles}). But which curvilinear curves? The curvilinear curves which are the perturbations of the two sides of the horn when perturbing $\eps$ slightly from $\eps_n$. In particular, if the formal invariant of $\tau^0$ is nonzero, then so is the case for the formal invariant of $\tau_{\eps_n}$. And this is the case for all $\exp\left(\frac{2\pi i p}{q}\right)$. We see that a lot of information is encoded in $\tau^0$, i.e. in $\psi^0$! This explains why the modulus is so large.

\section{The irregular nonresonant singular point of codimension 1 (Example 7)}\label{sec:Example7}

This case corresponds to the confluence of two regular singular points of a linear differential system. A  normal form for a generic unfolding is given by 
$$(x^2-\eps)\frac{dy}{dx} =
\left(D_0(\eps)+ D_1(\eps)x\right) y, \qquad y\in\C^n,$$
where $D_0(\eps)={\rm diag} (\lambda_1(\eps), \dots, \lambda_n(\eps))$ with the $\lambda_j(\eps)$ all distinct and  $D_1(\eps)={\rm diag} (\nu_1(\eps), \dots, \nu_n(\eps))$  is also diagonal. 
Modulo a rotation in $x$, we can always suppose that 
\begin{equation}{\rm Re}(\lambda_1(0))>\dots >{\rm Re} (\lambda_n(0)).\label{cond_order}\end{equation}

It can be shown (see \cite{LR}) that the system itself is analytically conjugate in a neighborhood of the origin to a system of the form
$$(x^2-\eps)\frac{dy}{dx} =
\left(D_0(\eps)+ D_1(\eps)x+ (x^2-\eps)B(x,\eps)\right) y, \qquad y\in\C^n.$$

For $\eps\neq0$, there are two regular singular points at $\pm \sqrt{\eps}$ with respective eigenvalues 
$$\mu_j^{\pm}(\eps)= \pm\frac{\lambda_j(\eps)\pm \nu_j(\eps)\sqrt{\eps}}{2\sqrt{\eps}}.$$
The monodromy map at a singular point is defined as follows: we consider a loop $\gamma: [0,1]\rightarrow \C$ surrounding the singular point, and an initial condition $y(0)\in \C^n$. Let $y(t)$ be the solution of the differential equation along the loop. Then the monodromy map around the singular point is the linear map $M$: $y(0)\mapsto y(1)$. Its similarity class is independent of the free homotopy class of $\gamma$.    The eigenvalues of the monodromy maps around each singular point are given by $\exp(2\pi i \mu_j^\pm(\eps))$. 

\bigskip

We can already observe the following:
\begin{enumerate} 

\item The eigenvalues of the monodromy have a very wild behaviour when $\eps\to 0$, namely that of an essential singularity in $\sqrt{\eps}$. But this wild behavior depends only on the formal normal form.
\item The eigenvalues of the monodromy map are distinct for generic values of the parameters, in which case the monodromy maps are diagonalizable: for these parameter values there exists at each singular point a basis of \emph{eigensolutions}, i.e. of special solutions, which are eigenvalues of the corresponding monodromy maps. This is the case where the regular singular points are \emph{nonresonant}. 
\item Of course, we expect that the generic situation is when the eigensolution at $-\sqrt{\eps}$ attached to the eigenvalue $\exp(2\pi i \mu_j^-(\eps))$ is not the analytic extension of the eigensolution at $\sqrt{\eps}$ attached to the eigenvalue $\exp(2\pi i \mu_j^+(\eps))$. 
\item Moreover, for special sequences of resonant values of the parameter converging to the origin, at least two eigenvalues of one monodromy map are equal, and the monodromy map is in general not diagonalizable.  Then  some solutions have logarithmic terms. 
\end{enumerate}

\bigskip 

In the same spirit as for the saddle-node we can go further:

\bigskip

\noindent{\bf Conclusion 1.} If we have divergence of the change of coordinate to the normal form when $\eps=0$, then this divergence 
in the limit forces the
eigenbases at each singular point to mismatch for sufficiently small $\eps$. 

\bigskip
\noindent{\bf Conclusion 2.} If we have divergence of the change of coordinate to the normal form when $\eps=0$, this will force the existence of sequences of resonant parameter values converging to $0$ for which there will exist solutions with nonzero logarithmic terms. This is again the  \emph{parametric resurgence phenomenon}. 

\bigskip

\noindent{\bf Conclusion 3.} In the formal normal form of the unfolding the number of parameters is equal to the number of analytic invariants at the linear level (here the eigenvalues) at each simple singular point.

\bigskip

Let us now discuss in more detail the case $\eps\neq0$. The eigenvalue $\mu_j^+(\eps)$ at $\sqrt{\eps}$ is almost opposite to the eigenvalue $\mu_j^-(\eps)$ at $-\sqrt{\eps}$. And the eigenvalues control the asymptotic behavior of the eigensolutions at the singular points; indeed, there exists (generically), in the neighborhood of each singular point, a basis of solutions  $\{\gamma_{j,\eps}^\pm\}$ with asymptotic expansion \begin{equation}\gamma_{j,\eps}^\pm(x) \simeq (x\mp \sqrt{\eps})^{\mu_j^\pm(\eps)} \left(c_j(\eps) \mathbf{e}_j+ O(x\mp\sqrt{\eps})\right),\label{asympt}\end{equation} where $c_j(\eps)\in \C^*$. 
\bigskip

\begin{figure}\begin{center} 
\subfigure[$\eps>0$]{\includegraphics[width=3.5cm]{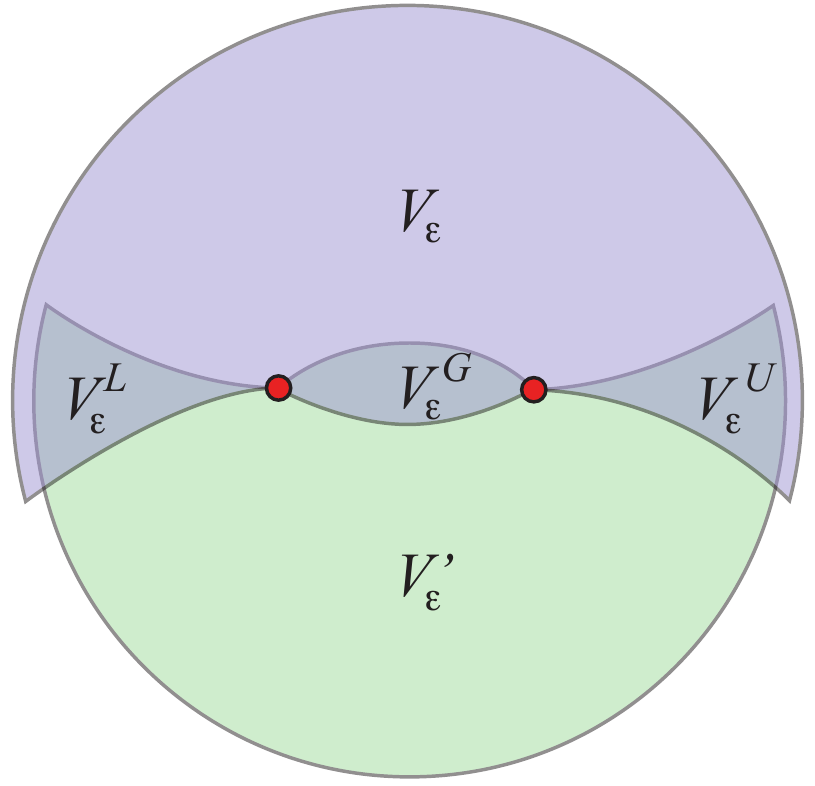}}\qquad \subfigure[$\eps=0$]{\includegraphics[width=3.5cm]{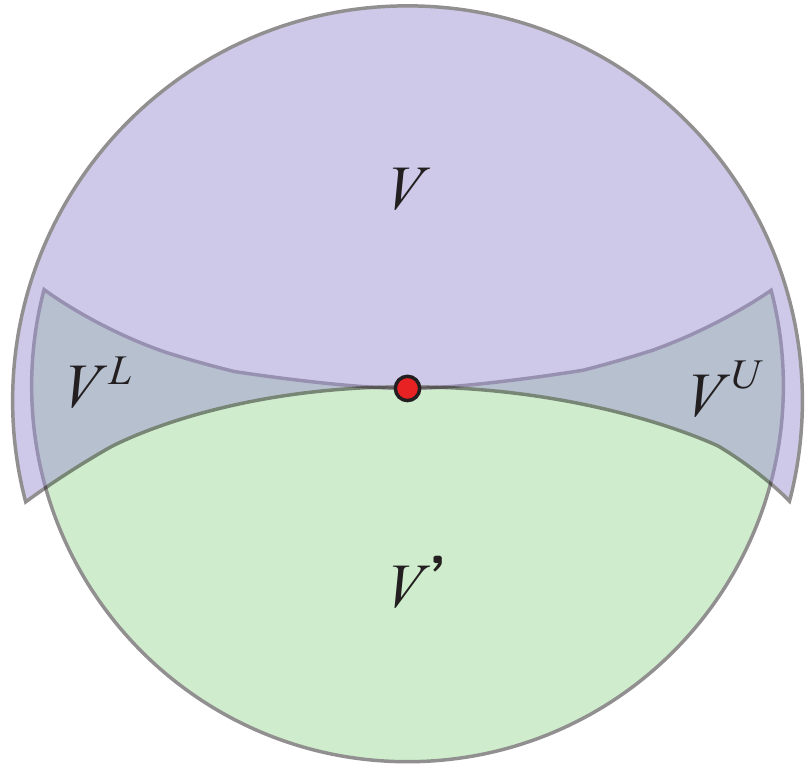}}\caption{The two sectors $V_\eps$ and $V_\eps'$ and their three (resp. two intersection parts) for $\eps>0$ (resp. $\eps=0$).}\label{fig:sectors} \end{center}\end{figure}

Let us take the case $\eps\in \R^+$. Then, using \eqref{cond_order}, this asymptotic behavior induces an ordering of the eigensolutions according to flatness
$$\begin{cases} 
\gamma_{1,\eps}^+\prec \gamma_{2,\eps}^+\prec \dots\prec \gamma_{n,\eps}^+,\\
\gamma_{n,\eps}^-\succ \dots\succ \gamma_{2,\eps}^-\succ \gamma_{1,\eps}^-.
\end{cases} 
$$
In turn, this allows to define flags of the solution  space at each singular point:
$$\begin{cases} 
W_{1,\eps}^+\subset W_{2,\eps}^+\subset \dots\subset W_{n,\eps}^+,\\
W_{n,\eps}^-\supset \dots\supset W_{2,\eps}^-\supset W_{1,\eps}^-.
\end{cases} 
$$ 

The flags depend analytically on $\eps$, with a continuous limit at $\eps=0$. From the nonresonance of the irregular singular point at $\eps=0$, it follows that the flags are transversal at $\eps=0$, yielding that  they are transversal for small $\eps>0$. Hence, taking two generalized sectors $V_\eps$ and $V_\eps'$ covering a disk $\D_r$ as in Figure~\ref{fig:sectors}(a), this allows to define on $V_\eps$ (resp. $V_\eps'$) a basis $\mathcal B_\eps$ (resp. $\mathcal B_\eps'$) that has the right asymptotic behavior at each singular point, namely the $j$-th vector has the asymptotic behavior \eqref{asympt} at $\pm \sqrt{\eps}$. 
Moreover, the bases $\mathcal B_\eps$ and $\mathcal B_\eps'$ are almost unique, the only degree of freedom being nonzero multiples. It is possible to choose bases depending analytically on $\eps\neq0$ in a sector around $\R^+$, and with continuous limit at $\eps=0$. 

 \bigskip
 
\begin{figure}\begin{center} \includegraphics[width=3.5cm]{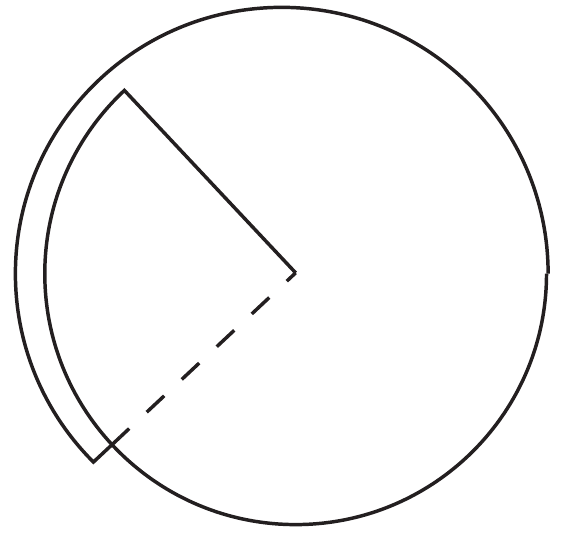}\caption{The sector $\Omega$ in $\hat{\eps}$-space.}\label{fig:epsilon}\end{center}\end{figure}

 Can we push these bases continuously for all values of $\eps$? The surprise is that the answer is positive, and we can construct bases $\mathcal B_{\hat{\eps}}$ (resp. $\mathcal B_{\hat{\eps}}'$) depending analytically on $\hat{\eps}$ in a sector of opening larger than $2\pi$ in the universal covering of $\eps$-space punctured at the origin (Figure~\ref{fig:epsilon}), and with a continuous limit at $\eps=0$. The key ingredient is the following lemma, which can be proved by mere calculations.
 
\begin{lemma} Let $\beta = a+ib\in \C$, with $b\neq0$. Then there exists logarithmic spirals $\gamma$ and $\gamma'$ with limit point at $0$ such that $\lim_{\begin{subarray}{l}z\to 0\\z\in \gamma\end{subarray}}z^\beta=0$ and $\lim_{\begin{subarray}{l}z\to 0\\z\in \gamma'\end{subarray}}z^\beta=\infty$.\end{lemma}

\begin{figure}\begin{center} \includegraphics[width=8cm]{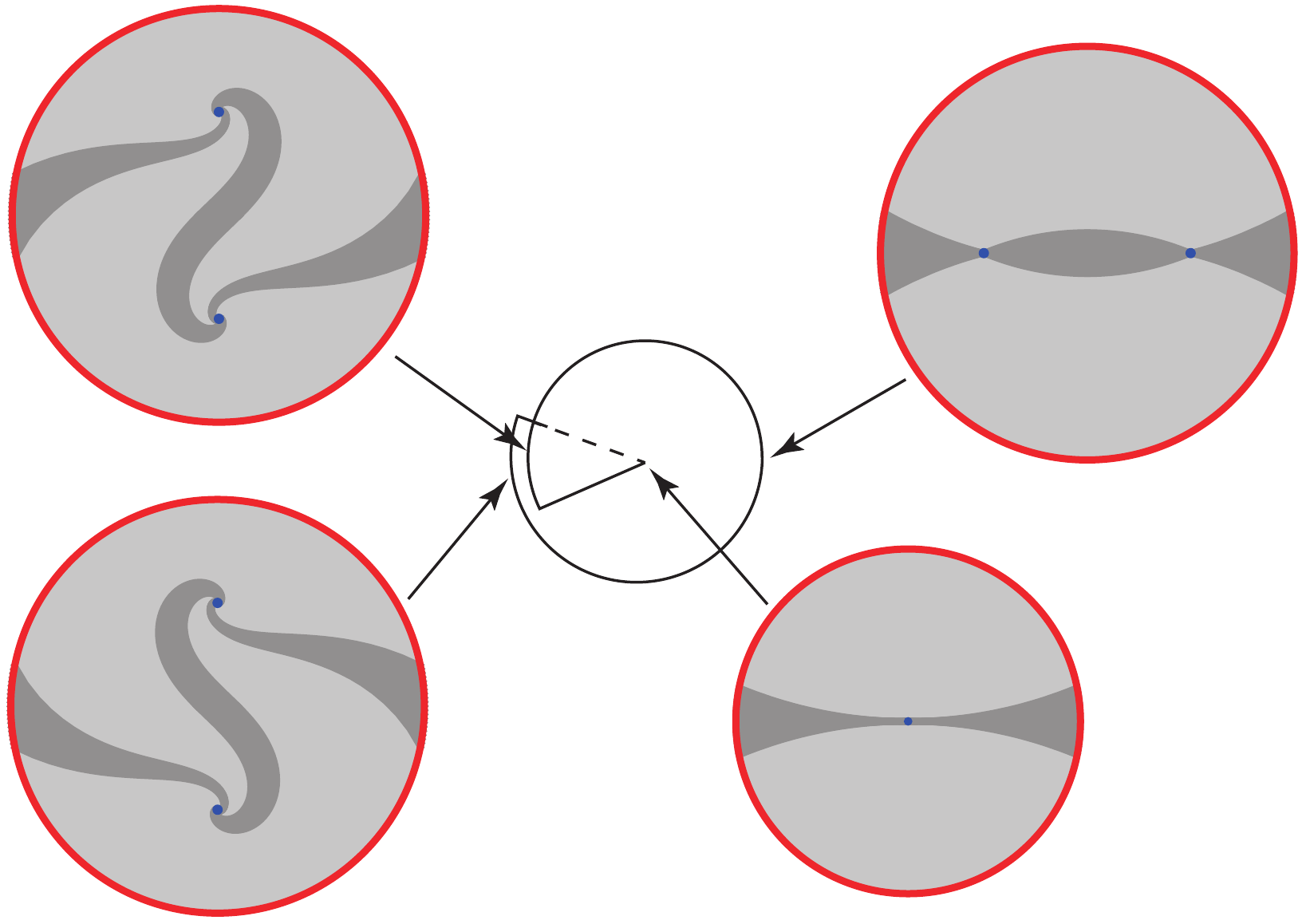}\caption{The two spiraling sectors for all values of $\eps$ and their three (resp. two) intersection parts in darker gray for $\eps\neq0$ (resp. $\eps=0$).}\label{fig:sectors_unfold} \end{center}\end{figure}

Then, when $\eps$ moves outside $\R^+$, we can deform the generalized sectors $V_\eps$ and $V_\eps'$ so that they approach the singular points $\pm\sqrt{\eps}$ along trajectories of some rotated vector field $\dot z = e^{i\alpha} (z^2-\eps)$, which are very close to logarithmic spirals (see Figure~\ref{fig:sectors_unfold}). 

Hence, if the sectors $V_\eps^U$,  $V_\eps^L$ and $V_\eps^G$ are the three intersection parts of the sectors $V_\eps$ and $V_\eps'$, then the comparison between the bases  $\mathcal B_\eps$ and $\mathcal B_\eps'$  over $V_\eps^U$,  $V_\eps^L$ and $V_\eps^G$
is the \emph{classifying object}. Because of the flatness properties, comparing the two bases is done  via an upper triangular matrix $S_\eps^U$ over $V_\eps^U$, a lower triangular matrix $S_\eps^L$ over $V_\eps^L$, and a diagonal matrix $S_\eps^G$ over $V_\eps^G$. 

\bigskip
\noindent{\bf Conclusion 5.} We have decomposed the dynamics into: \begin{itemize} 
\item a wild diagonal part $S_\eps^G$ depending only on the formal normal form and having no limit at $\eps=0$;
\item two triangular parts $S_\eps^U$ and $S_\eps^L$, which converge to the classical \emph{Stokes matrices}  $S_0^U$ and $S_0^L$ at the irregular singular point for $\eps=0$ (see Figure~\ref{fig:sectors}(b)).\end{itemize}

\bigskip 
\begin{remark} Since the change of basis is given by a diagonal matrix $S_\eps^G$ over $V_\eps^G$, we could have used the same basis $\mathcal B_\eps$ on both $V_\eps$ and $V_\eps'$, and this basis would have not been defined in a uniform way when turning around each singular point. In that point of view the dynamics is defined by the monodromy around the singular points. The price to pay is that this monodromy has no limit when $\eps\to 0$: its eigenvalues at each singular point are essential singularities in $\sqrt{\eps}$. \end{remark}

 \bigskip

\noindent{\bf Conclusion 6.} The construction can be pushed to cover all values of $\eps$  in a small neighborhood of the origin
 (see Figure~\ref{fig:sectors_unfold}), with the following constraints:
\begin{itemize}
\item The generalized sectors $V_\eps$ and $V_\eps'$ may have to spiral at the singular points;
\item The construction is ramified in $\eps$, i.e. done on a sector $\Omega$ of opening larger than $2\pi$ centered on $\R^+$ as in Figure~\ref{fig:epsilon}.
\end{itemize}

\bigskip

As before, let us denote by $\hat{\eps}$ an element of the universal covering of the $\eps$-space punctured at $0$.  

\begin{theorem}\label{thm:LR} \cite{LR} Two families of linear differential systems unfolding an irregular nonresonant singular point of Poincar\'e rank $1$ are analytically conjugate if and only if they have the same \emph{modulus}, namely
\begin{enumerate} 
\item they have the same formal normal form;
\item they have \lq\lq equivalent\rq\rq\  collections of unfolded Stokes matrices $(S^U_{\hat{\eps}},S^L_{\hat{\eps}})$ depending continuously on  $\hat{\eps}$ in a sector  $\Omega =\{\hat{\eps}\; : \: |\eps|<\rho; \arg(\hat{\eps})\in(-\pi -\delta, \pi +\delta)\}$ for some $\delta\in (0,\pi)$ (see Figure~\ref{fig:epsilon}). The equivalence on collections of unfolded Stokes matrices $(S^U_{\hat{\eps}},S^L_{\hat{\eps}})$ is defined as follows:
$\left(S^U_{1,\hat{\eps}},S^L_{1,\hat{\eps}}\right)\equiv \left(S^U_{2,\hat{\eps}},S^L_{1,\hat{\eps}}\right)$ if and only if there exists invertible diagonal matrices $D(\hat{\eps}), D'(\hat{\eps})$ depending continuously on $\hat{\eps}\in \Omega$ such that
$$ \begin{cases} 
 S^U_{1,\hat{\eps}}= D(\hat{\eps})\; S^U_{2,\hat{\eps}}\;D'(\hat{\eps}),\\
 S^L_{1,\hat{\eps}}= D(\hat{\eps})\; S^L_{2,\hat{\eps}}\;D'(\hat{\eps}).\end{cases}$$ 
 \end{enumerate} 
\end{theorem}

\noindent {\bf Conclusion 7.}  The conditions for analytic conjugacy of two linear differential systems with an irregular singular point of Poincar\'e rank $k$ are well known in the literature: \emph{the two systems must have the same normal form and equivalent collections of Stokes matrices.} Theorem~\ref{thm:LR} states that the unfolding of the modulus  of the system for $\eps=0$ is the modulus of the unfolded system. 

\bigskip Again, conclusions similar to Conclusions 5-7  are a general feature for all our examples, as well as classification theorems of the type of Theorem~\ref{thm:LR}. 

\section{The common features to all examples when the codimension is $1$} 

Let us summarize  some common features to all generic unfoldings of the singularities of all examples when $k=1$.
\begin{enumerate} 
\item There exists a formal normal form, whose number of parameters is equal to the number of analytic invariants at the singular points, when these are simple.
\item Except for Example 7, the parameter of the formal normal form is an analytic invariant. See for instance \eqref{canonical_parameter} in the case of the saddle-node.
\item The description of the dynamics is not uniform in the parameter space. In Examples 1, 3, 6 (in the case of a zero angle) and 7, it is done over a sector $\Omega$ of opening $2\pi + 2\delta$ for some $\delta \in (0,\pi)$, in the universal covering $\hat{\eps}$ of the $\eps$-space punctured at $0$ (see Figure~\ref{fig:epsilon}). The upper bound for $\delta$ depends on the singularity type. It is equal to $\pi$ for parabolic points of diffeomorphisms, and saddles or saddle-nodes of vector fields. It is generically smaller for nonresonant singular points of linear differential systems, and depends on the eigenvalues $\lambda_j$. The larger $\delta$, the smaller the radius of the sector in $\hat{\eps}$. In Examples 2, 4, 5, 6 (in the case of a nonzero angle), the description of the dynamics is done over two sectors $\Omega_1$ and $\Omega_2$ of opening larger than $\pi$ and covering a disk in parameter space $\eps$. 
\item The dynamics of the generic $1$-parameter family of systems  is described on the union of  two generalized sectors $V_{\hat{\eps}}$ and $V_{\hat{\eps}}'$ covering a disk $\D_r$. The sectors are bounded by $C^0$ curves, each being a union of a finite number of flow lines of some rotated vector field $\dot z = e^{i\alpha}(z^2-\eps)$ of the \emph{organizing vector field} $\dot z = z^2-\eps$. 
 In Examples 2, 4, 5, 6 (in the case of a nonzero angle), we rather use the vector field $\dot z = z(z-\eps)$, which has a fixed singular point at the origin. 
\item There exists over each sector $V_{\hat{\eps}}$ and $V_{\hat{\eps}}'$ an almost unique change of coordinate to the formal normal form. 
\item Then the modulus is obtained by comparing the changes of coordinates to the normal form over the three parts $V_{\hat{\eps}}^U$,  $V_{\hat{\eps}}^L$ and $V_{\hat{\eps}}^G$ of the intersection of $V_{\hat{\eps}}$ and $V_{\hat{\eps}}'$ (see Figures~\ref{fig:sectors} and \ref{fig:sectors_unfold}). 
\item The changes of coordinates to the formal normal form over $V_{\hat{\eps}}$ and $V_{\hat{\eps}}'$ can be chosen so that their comparisons over $V_{\hat{\eps}}^U$ and  over $V_{\hat{\eps}}^L$ have a limit when ${\hat{\eps}}\to 0$, which is given by the classical modulus for the case $\eps=0$. Then the change of coordinate over $V_{\hat{\eps}}^G$ is trivial (diagonal or linear depending on the context), but with very wild behavior and no limit when ${\hat{\eps}}\to 0$: the nonzero entries have an  essential singularity in $\sqrt{\hat{\eps}}$ (resp. $\eps$) in Examples 1, 3, 6 (in the case of a zero angle) and 7 at $\eps=0$ (resp. in Examples 2, 4, 5, 6 (in the case of a nonzero angle)). However, this wild behavior is completely controlled by the formal normal form. 
\item Hence, we have a decomposition of the dynamics  into a wild formal part with no limit at $\eps=0$, and an analytic part, which tends to the analytic part of the modulus when $\hat{\eps}=0$.
\end{enumerate}

\subsection{The realization in codimension $1$} All the  classification problems for $\eps=0$ are classical (see for instance \cite{IY} or \cite{Z}). In all of them a complete modulus of classification is given, which has two parts: 
\begin{itemize}
\item a formal part depending on a finite number of parameters;
\item an analytic part, usually given as an equivalence class (modulo the global symmetries of the formal normal form).\end{itemize}  Moreover, the realization problem is solved, which consists in identifying which moduli can be realized. This allows giving the moduli set. 
Except for Example 7, this moduli set is of infinite dimension and no topology is given on it. We rather see theorems of the form: \emph{If the codimension $1$ (resp. $k$) singularity occurs in a family with fixed formal normal form depending analytically on $\ell$ parameters, then there exists a representative of the modulus depending analytically on these parameters.}

\bigskip
The realization is trickier when one is considering generic one-parameter unfoldings of a codimension $1$ singularity of the types described above. Remember that the analytic part of the modulus does not depend analytically on the parameter in a full neighborhood of the origin. It is either defined on a sector $\Omega$ of opening larger than $2\pi$ in the universal covering $\hat{\eps}$ of the parameter space punctured at the origin, or in two sectors $\Omega_j$, $j=1,2$, of opening larger than $\pi$. It is possible to realize any potential  modulus for each $\hat{\eps}$ and to get a family depending analytically on $\hat{\eps}\in \Omega$ (resp. $\hat{\eps}\in \Omega_j$, $j=1,2$),  with continuous limit at $\eps=0$. But, without an additional condition, there is no reason why there would exist a change of coordinates over $\Omega$ (resp. changes of coordinates over $\Omega_1$ and $\Omega_2$) depending analytically on $\hat{\eps}$ and transforming the realized family into a family depending analytically on $\eps$. An obvious necessary condition is that on the self-intersection of the neighborhood $\Omega$ in $\eps$-space (resp. on $\Omega_1\cap \Omega_2$) the two realizations are conjugate. Together with a technical condition on the limits when $\eps=0$ this turns out to be sufficient, and realization theorems exist in codimension $1$ for the seven types of singularities listed above (see for instance  \cite{CR} for the parabolic case, \cite{Ro3} for the fixed point with periodic multiplier and resonant saddle, and \cite{LR} for the irregular nonresonant singular point of a linear differential equation).

\section{Moving to higher codimension}

When we have the confluence of $k+1$ \emph{special objects} we often observe $k$-summability of the normalizing changes of coordinates. 
This is the case in our seven examples. 

\bigskip

The dynamics of the vector field $\dot z=z^2-\eps$ (or $\dot z = z(z-\eps)$) played a very important role in the codimension $1$ case. It the codimension $k$ case, it is replaced by the dynamics of the vector field
$\dot z = P_\eps(z)$, with $P_\eps$ given in \eqref{def:P}, or by that of $\dot z = zQ_\eps(z)$, with $Q_\eps$ given in \eqref{def:Q}. Indeed, in Examples 1-6, there is an underlying $1$-dimensional map, whose formal normal form is the time-one map of a vector field close to $\dot z = P_\eps(z)$ or $\dot z = zQ_\eps(z)$ for small $(z,\eps)$. We limit the discussion to the case of $\dot z = P_\eps(z)$.

\bigskip
In Example 7, a generic unfolding has the form 
$$P_\eps(x)\frac{dy}{dx} = A_\eps(x) y, \qquad y\in \C^n,$$ wich can also be rewritten as an ODE
\begin{align} \begin{split} 
\dot x &= P_\eps(x),\\
\dot y &=A_\eps(x) y,\end{split}\label{LS}\end{align}
in which the vector field $\dot x = P_\eps(x)$ organizes the type of the singular points. 
\bigskip

\begin{figure}\begin{center} 
\subfigure[The pole at infinity]{\includegraphics[width=3.5cm]{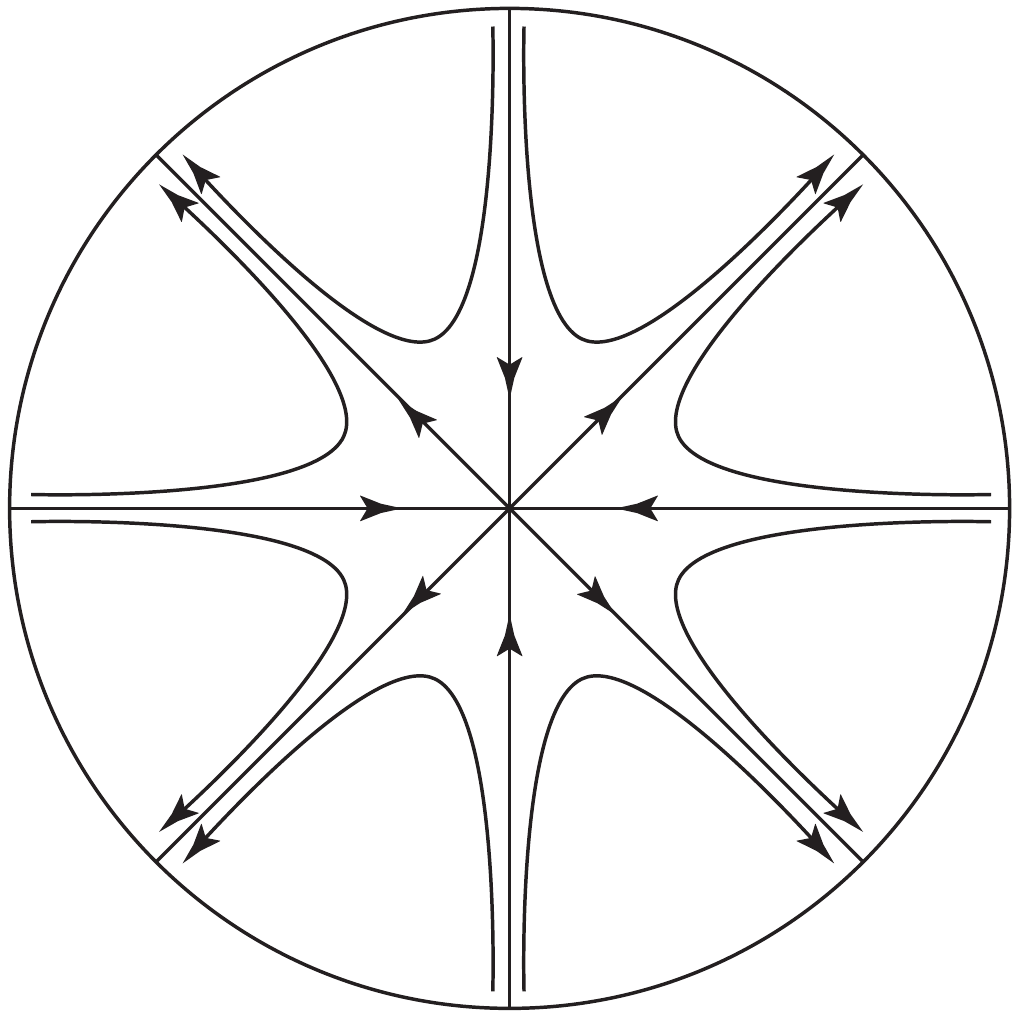}}\quad\subfigure[The separatrix graph]{\includegraphics[width=3.5cm]{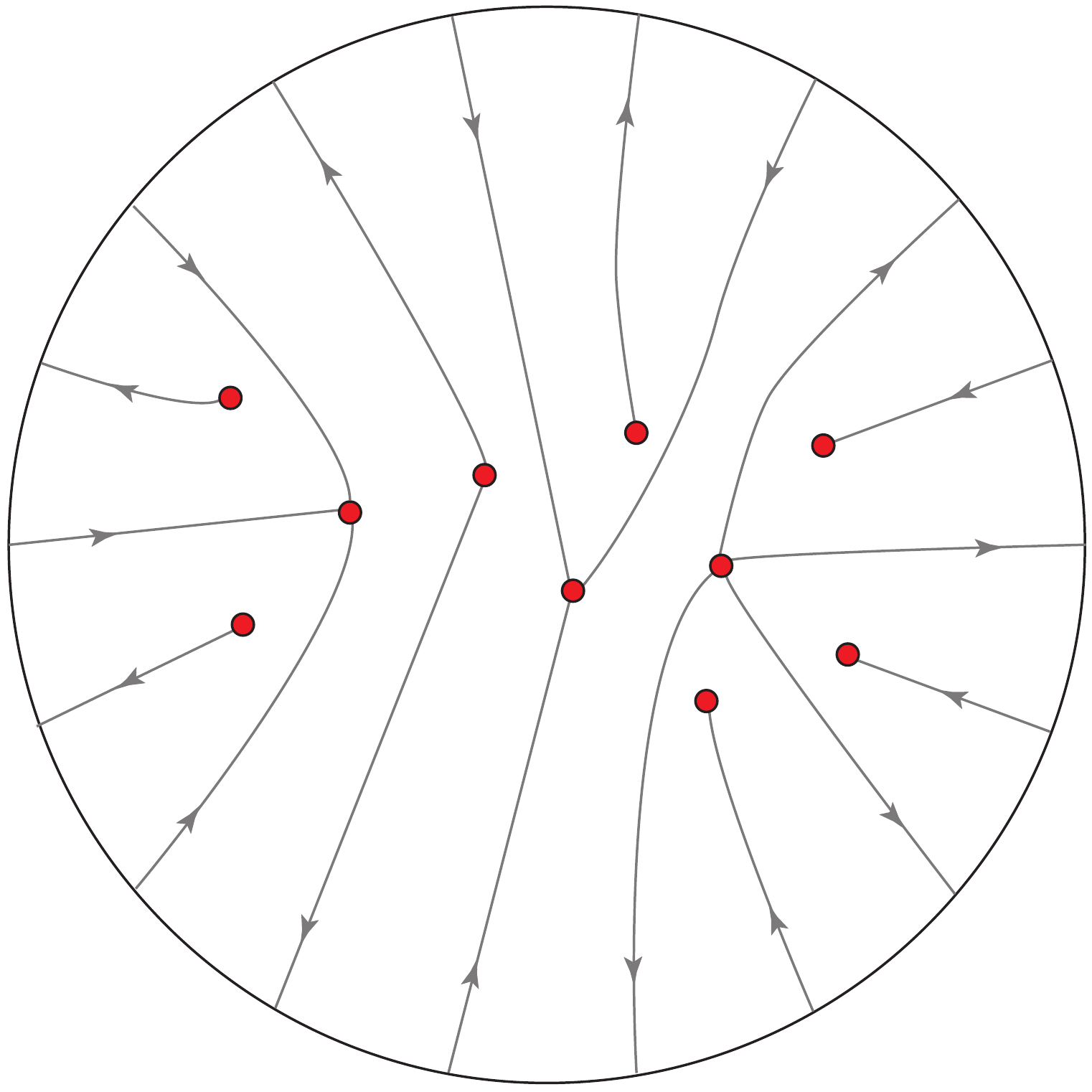}}\quad\subfigure[Two zones]{\includegraphics[width=3.5cm]{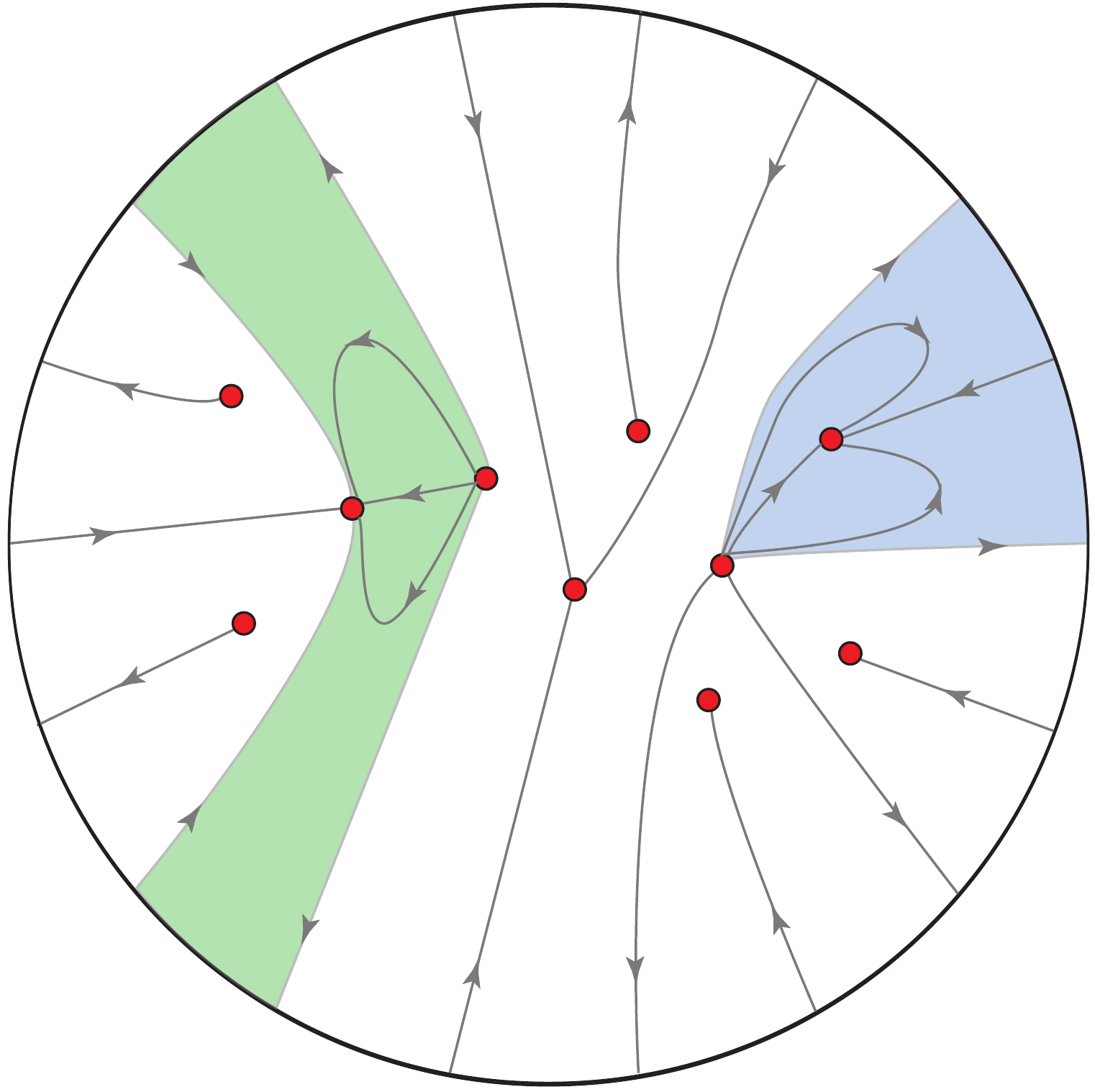}}
\caption{The phase portrait of $\dot z = P_\eps(z)$: (a) the pole at infinity for $k=4$ and the $8$ separatrices; (b) a separatrix graph and (c) two zones for $k=9$.}\label{phase_portrait}\end{center}\end{figure}

The  dynamics of the family of vector fields $\dot z = P_\eps(z)$ was studied by Douady, Estrada and Sentenac in a visionary paper \cite{DES}.
For each parameter value, the dynamics  on $\CP^1$ is governed by the pole of order $k-1$ at infinity, which has $2k$ separatrices, in turns attracting or repelling (Figure~\ref{phase_portrait}(a)). Generically, these separatrices land at the singular points (Figure~\ref{phase_portrait}(b)). Exceptionally, a homoclinic loop occurs between two separatrices. The structurally stable vector fields are dense. These are the vector fields  with simple singular points and no homoclinic loop through infinity. The structurally stable vector fields are of $C(k)$ different topological types,  where $C(k) = \frac{\binom{2k}{k}}{k+1}$ is the $k$-th Catalan number. The parameter values corresponding to a given topological type is a simply connected open domain.
For parameter values in any of these domains, the \emph{separatrix graph} composed of the union of the separatrices divides $\CP^1$ into $k$ simply connected domains in $z$-space called \emph{zones} (Figure~\ref{phase_portrait}(c)). Each zone has no singular point inside and two singular points on the boundary, both  of node or focus type, one attracting, one repelling. Moreover,  all trajectories inside a zone have their $\alpha$-limit at the repelling singular point, and their $\omega$-limit at the attrating singular point.  
\bigskip

In the codimension $1$ case, there was no uniform way to describe the dynamics for all parameter values: the dynamics needed to be described over a sector of opening larger than $2\pi$ in the universal covering of $\eps$-space at the origin. In higher codimension, the dynamics will be described over $C(k)$ open sets $U_j$, called DES-domains, which are enlargements of the parameter regions described by Douady, Estrada and Sentenac, so that the union of the $U_j$ cover all parameter values where the special objects (singular points, periodic orbits, etc.) are distinct. The open sets $U_j$ are roughly unions of domains of local structural stability of a given topological type for rotated vector fields $\dot z = e^{i\alpha}P_\eps(z)$.
\bigskip

\begin{figure}\begin{center}\includegraphics[width=4cm]{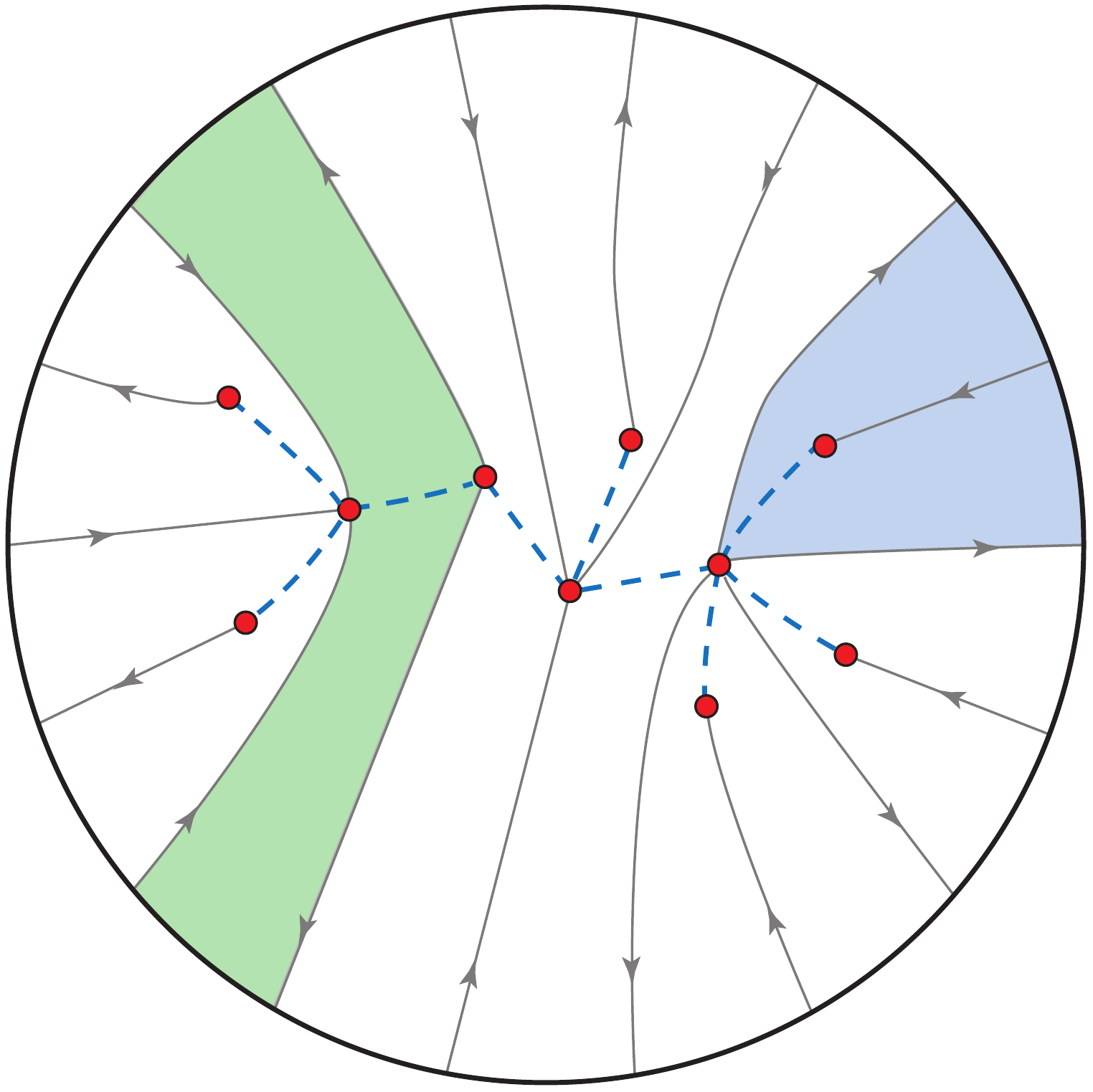}\caption{The zones are cut into half-zones along the dotted lines.}\label{half-zones}\end{center}\end{figure}

This description will allow to generalize to codimension $k>1$ what has been done in codimension $1$. Indeed, on each zone we will look for almost unique changes of coordinates to normal form, very similar to the ones obtained in the codimension $1$ case. Comparing these normalizing changes of coordinates will yield the modulus of analytic classification, which encodes the obstruction to a global analytic change of coordinate to normal form. In practice, if we do this however, each zone splits into two parts at the limit when the two attached singular points merge together, and if we work with a unique change of coordinate to the normal form we cannot have convergence at the limit when the two points merge together. Hence, we will work with half-zones as in Figure~\ref{half-zones}. 
The change of normalizing coordinate from one half-zone to the other one along the border line is essentially trivial: it is a symmetry of the formal normal form over the full neighborhood and it depends only on the formal normal form. Let us now give more detail in the special case of Example 7. 
\bigskip

\subsection{Example 7 revisited in codimension $k$.}  Let us consider a nonresonant irregular singular point of Poincar\'e rank $k$  of a linear differential system

$$x^{k+1}\frac{dy}{dx} = A(x) y, \qquad y\in \C^n.$$
The  formal normal form of the unfolding is given by $$P_\eps(x)\frac{dy}{dx} = (D_0+ D_1x+\dots D_kx^k)
y$$ where
$P_\eps$ is given in \eqref{def:P} and the $D_j$ are diagonal. We can of course suppose that $A(0)=D_0$ is diagonal and that the eigenvalues satisfy \eqref{cond_order}.

All together this gives $(k+1)n$ parameters to control the eigenvalues at  $k+1$ singular points.

\bigskip

For $k=1$, the almost unique bases we had constructed came from flags of solutions associated to each singular point. The flags were in inverse direction and transverse one to the other. 
When $k>1$, let us first consider the case where all singular points of the vector field $\dot z =P_\eps(z)$ have real eigenvalues and are hence attracting or repelling nodes. In that case we work over each zone, which is attached to one attracting singular point $z_j$ and one repelling singular point $z_\ell$ of the vector field. On each zone we find one flag attached to each singular point. The flags are inverse one to the other and transversal. Hence, by intersecting the flags associated to the singular points $z_j$ and $z_\ell$  we can,  as in the codimension $1$ case, find on each zone a basis of solutions $\mathcal{B}_\eps=\{\gamma_{1,\eps},\dots, \gamma_{n,\eps}\}$ such that 
\begin{equation}\begin{cases} 
\gamma_{1,\eps}\prec \gamma_{2,\eps}\prec \dots\prec \gamma_{n,\eps},& \text{near}\:z_j,\\
\gamma_{n,\eps}\succ \dots\succ \gamma_{2,\eps}\succ \gamma_{1,\eps},& \text{near}\:z_\ell.
\end{cases}\label{asymp_order}
\end{equation}
Moreover this basis is unique up to multiples of the basis vectors.

When we move to values of $\eps$ for which the eigenvalues of singular points of $\dot z =P_\eps(z)$  are no more real we may have to replace the zones by \emph{generalized zones}. In a generalized zone we may approach $z_j$ along trajectories of some $\dot z = e^{i\alpha_j} P_\eps(z)$ and approach $z_\ell$ along trajectories of $\dot z = e^{i\alpha_\ell} P_\eps(z)$, where $\alpha_\ell $ could possibly be different from $\alpha_j$. When approaching along such trajectories which are very close to logarithmic spirals we deform continuously the flags while preserving their inverse character and their transversality. We can then deform the basis  $\mathcal{B}_\eps$ associated to a (generalized) zone, while preserving \eqref{asymp_order}. Each angle $\alpha_u$ is attached to its singular point $z_u$, and all generalized zones adherent to a fixed $z_u$ use the same angle.  The comparison of the bases over different zones are given by the \lq\lq Stokes matrices\rq\rq, which are unfoldings of the classical Stokes matrices. 

\begin{theorem} \cite{HLR} Two families unfolding an irregular nonresonant singular point of Poincar\'e rank $k$ are analytically conjugate if and only if 
\begin{enumerate} 
\item They have the same formal normal form;
\item Over each DES domain $U_s$ they have  \lq\lq equivalent\rq\rq\  collections of unfolded Stokes matrices $\mathcal{S}_s=\left(S^U_{1,s, \eps},\dots, S^U_{k,s,\eps},S^L_{1,s,\eps},\dots,S^L_{k,s,\eps}\right)$. The   $C(k)$ open DES domains $U_s$ cover the parameter space minus the discriminant set. The equivalence relation on collections $\mathcal{S}_s$ of Stokes matrices over $U_s$ is that corresponding to the degrees of freedom in the choice of bases over the different zones. \end{enumerate}\end{theorem}

\bigskip

Moreover, given a formal normal form,  the moduli space, i.e. the realizable collections of unfolded Stokes matrices $\left\{\left\{\left(S^U_{1, j,\eps},\dots, S^U_{k, j,\eps},S^L_{1,j,\eps},\dots,S^L_{k,j,\eps}\right)\right\}_{\eps\in U_j}\right\}_{j=1}^{C(k)}$ has been identified \cite{HR}. The realization is done in the following way. One first realizes the modulus over each DES domain $U_j$ in parameter space. A necessary condition, the \emph{compatibility condition} is then introduced which guarantees that two realizations over DES domains $U_j$ and $U_{j'}$ are analytically equivalent over $U_j\cap U_{j'}$: the condition is that the monodromy representations given by the formal normal form and the unfolded Stokes matrices are conjugate, plus a technical condition: the cocycle giving the conjugacy is trivial. Together with limit conditions when approching the discriminantal set, these conditions are also sufficient. 

\subsection{The first six examples.}Classification theorems of the type \cite{HLR} exist or are straightforward for all seven examples. In all cases, the modulus of analytic classification is given by the formal normal form and $C(k)$ unfoldings of the modulus over $C(k)$ DES domains in parameter space. One difficulty was to identify canonical parameters (see Section~\ref{sec:can_par} below). 
  The realization is trickier and only done in Example 7. A partial realization is done in the saddle-node case, under the hypothesis that the unfolding has an analytic center manifold for all $\eps$ (see \cite{RT2}). It is easy to realize unfolded moduli over each DES domain. The difficulty is to identify the compatibility conditions that the different unfolded moduli over the $C(k)$ DES domains in parameter space must satisfy in order that there exists a realization which is analytic in the parameters. For the saddle-node of codimension $k$, in the particular case of an analytic center manifold, there exists a monodromy pseudo-group and the compatibility condition is similar to that of Example 7. Indeed, Example 7 can be rewritten in the form \eqref{LS}, which has some similarity with the unfolding of a saddle-node with analytic center manifold. 

\subsection{Canonical parameters.}\label{sec:can_par}

Showing that there exist canonical parameters for Examples 1-6 in the higher codimension case is non trivial. In all cases, it amounts to show that the unfolding 
$\dot z = \frac{P_\eps(z)}{1+a(\eps)z^k}$ ($P_\eps$ defined in \eqref{def:P}) of a parabolic singular point $\dot z = z^{k+1}+ o\left(z^{k+1}\right)$ is universal. 
The versality was proved by Kostov (\cite{Ko}). The universality was first proved in \cite{RT1}. The precise statement is the following. 

\begin{theorem}
\label{parametres_inv} Let  $\Psi\,:\,\left(x,\varepsilon\right)=\left(x,\varepsilon_{0},\ldots,\varepsilon_{k-1}\right)\mapsto\left(\varphi_{\varepsilon}\left(x\right),h_{0}\left(\varepsilon\right),\ldots,h_{k-1}\left(\varepsilon\right)\right)=\left(z,h\right)$ be a germ of an analytic change of
coordinates
at $\left(0,0,\cdots,0\right)\in\mathbb{C}^{1+k}$. The following
assertions are equivalent: 
\begin{enumerate}
\item The families $X_\eps: \left(\frac{P_{\varepsilon}\left(x\right)}{1+a\left(\varepsilon\right)x^{k}}\frac{\partial}{\partial x}\right)_{\varepsilon}$
and $\widetilde{X}_\eps: \left(\frac{P_{h}\left(z\right)}{1+\tilde{a}\left(h\right)z^{k}}\frac{\partial}{\partial z}\right)_{h}$
are conjugate under $\Psi$. 
\item There exist $T\in\C\{\varepsilon\}$, and $\tau$ with $\tau^{k}=1$,
such that

\begin{itemize}
\item $\varphi_{\varepsilon}\left(x\right)=\Phi_{X_{\eps}}^{T(\eps)}\circ R_{\lambda}\left(x\right)$,
where $R_{\tau}(x)=\tau x$ and $\Phi_{X_{\eps}}^{T(\eps)}$ is the flow of $X_\eps$ at time $T(\eps)$, 
\item $\varepsilon_{j}=\tau^{j-1}h_{j}\left(\varepsilon\right)$, 
\item $a\left(\varepsilon\right)=\tilde{a}\left(h\left(\varepsilon\right)\right)$. 
\end{itemize}
\end{enumerate}
\end{theorem}
The proof is done by infinite descent (see  \cite{RT1}, Theorem 3.5, or \cite{KR}).

\section{Perspectives}

In the case of dynamical systems having simple normal forms, we encounter different types of divergence of the normalizing changes of coordinates, for instance: 
\begin{enumerate}
\item {\bf $k$-summability.} It usually occurs in $1$-resonant systems, i.e. all resonance relations are generated by one rational relation between the eigenvalues.
\item {\bf Multi-summability.} It occurs when studying irregular resonant singular points of  linear differential equations.
\item{\bf Small divisors.} It occurs close to an infinite number of independent resonances.\end{enumerate}

The divergence of the normalizing changes of coordinates reflects that the \emph{geometry} of the system is more complex than that of the normal form. In the case of $1$-resonance presented in this paper this comes from the coallescence of $k+1$ special \lq\lq objects\rq\rq. When unfolding to separate the objects, we can generically analytically transform to the normal form in the neighborhood of each special object and the divergence comes from the mismatch of these local models. 

The challenge is to find similar descriptions for the other cases and then to integrate all cases in a larger portrait. Indeed, as illustrated in this paper, unfolding a singularity allows to understand the local geometry of a system and the obstructions to the convergence of the normalizing transformation. 

The third case was already discussed above in Remark~\ref{rem:Yoccoz}. A fixed point of a $1$-dimensional germ of diffeomorphism with multiplier $e^{2\pi i \alpha}$ for irrational $\alpha$ is formally linearizable. It is analytically linearizable if $\alpha$ is Diophantian, and generically non linearizable if $\alpha$ is Liouvillian. Yoccoz showed that the precise frontier is the Bruno condition on the continued fraction of $\alpha$ (see \cite{Y}). He also showed that when we are very close to the Bruno condition, there is only one mechanism which prevents linearizability, namely the existence of an infinite number of periodic points. Generically the periodic points of period $q$ merge with the fixed point when $\alpha$ is perturbed to some close rational $\frac{p}{q}$. When $\alpha$ is more Liouvillian (closer to the rationals) P\'erez-Marco (\cite{PM1} and \cite{PM2}) showed  that there exist other mechanisms preventing orbital linearizability: the existence of a complicated invariant set, the \emph{hedgehog}. The geometric studies of (1) and (3) have been done in parallel, but would deserve more integration. 

In the case of an irregular singular point of Poincar\'e rank $k$  of a linear differential system, we have $k$-summability when the singular point is nonresonant. In the resonance case however we have multi-summability of the normalizing changes of coordinates. This is a more complex case, which does not share the common features described in this paper. In particular,  the dynamics of the vector field $\dot x = P_\eps(x)$ is no more organizing the geometry. A very interesting question is to study the geometric obstructions to convergence in the unfoldings of some system for which the normalizing transformation is multi-summable. Martin Klime\v{s} (see \cite{Ki}) determined the modulus of an unfolding of a resonant irregular singular point of Poincar\'e rank $1$ of a linear differential system of dimension $2$, but it is still premature to draw common features for the resonance case.

\end{document}